%% file: main.tex
\documentclass{siamart251216}
\usepackage{amssymb,mathtools}
\usepackage[nopatch=eqnum]{microtype}
\usepackage{bm}
\usepackage{braket}
\usepackage{faktor}

\newcommand{\R}{\mathbb{R}}
\newcommand{\N}{\mathbb{N}}
\newcommand{\B}{\mathcal{B}}

\newcommand{\HH}{\mathcal{H}}

\newcommand{\K}{\mathcal{K}}
\newcommand{\Pc}{\mathcal{P}}
\newcommand{\X}{\mathcal{X}}
\newcommand{\Ext}{{\rm Ext}}
\newcommand{\dd}{\mathrm{d}}
\newcommand{\qandq}{\quad \text{and} \quad }

\newcommand{\hilbert}{\mathcal{H}}
\newcommand{\purest}{\mathbf{P}_\hilbert}
\newcommand{\dens}{\mathfrak{S}_{\hilbert}}

\newcommand{\psDistance}{pseudo-distance }
\newcommand{\semiDistance}{semi-distance }
\newcommand{\semiMetric}{cost }
\newcommand{\pRoof}[1]{d_{*#1}}

\renewcommand{\leq}{\leqslant}
\renewcommand{\geq}{\geqslant}

\newsiamremark{remark}{Remark}

\title{Metric convex extensions and optimal transport on classical and quantum structures\thanks{Previous versions of this document were entitled ``\emph{Folded optimal transport and its applications to separable quantum optimal transport}''.\funding{This work was funded by the ERC Starting Grant project HighLEAP number 101077204 awarded to Virginie Ehrlacher.}}}
\author{Thomas Borsoni\thanks{CERMICS, \'Ecole des Ponts, Institut Polytechnique de Paris, and Inria, France (\email{thomas.borsoni@enpc.fr}).}}
\headers{Metric convex extensions and optimal transport}{T. Borsoni}

\begin{document}

\maketitle

\begin{abstract}
Convex roofs are natural ways to extend a scalar function defined on the extreme points of a convex set $K$ to the whole set. In this work we focus on the case $K = C \times C$ for some convex set $C$, and the extension of a distance $d$ defined on the extreme points of $C$. In general, the triangle inequality may not survive the extension process. We therefore consider the metric envelope of such extensions and study the metric properties they endow $C$ with. More generally, we also consider products of different convex sets $C_1 \times C_2$ and costs defined on their extreme points, in which case no metric envelope is involved. We then show a two-way correspondence between convex extensions and Kantorovich optimal transport formulations. On the one hand, these extensions can be recovered as quotients of standard, classical Kantorovich costs or Wasserstein distances by the barycenter map. On the other hand, various Kantorovich formulations re-write as such convex extensions. Classical optimal transport is recovered as the case of the simplex, and both the semiclassical optimal transport cost of Golse and Paul and the entanglement-free quantum Wasserstein objects of Beatty and Stilck Fran\c{c}a are recovered by considering the non-simplex convex set of quantum states.

This work provides a clear bridge between the theory of convex extensions and that of optimal transport, with a unifying framework for classical, semiclassical and entanglement-free quantum instances of Kantorovich formulations.
\end{abstract}

\begin{keywords}
convex extension, convex roof, Choquet theory, metric envelope, Wasserstein distance, quantum optimal transport
\end{keywords}

\begin{AMS}
46A55, 49Q22, 81P16, 52A07
\end{AMS}

\input{introduction}

\input{Section_convex_roofs}

\input{Section_proofs}

\input{Section_applications}

\section*{Acknowledgments}
I most warmly thank Virginie Ehrlacher and Genevi\`eve Dusson for our discussions from the beginning of the project, support, careful reading and many comments that helped improve this article. I kindly thank G\'eza T{\'o}th, J\'ozsef Pitrik and Pierre-Cyril Aubin for pointing out relevant references. I greatly thank an anonymous referee for the Journal of Statistical Physics for the significant and detailed remarks which initiated my re-work on this article. I deeply thank Vincent Boulard for our numerous conversations and his insights, which greatly helped me reshape the article and find numerous applications of the present work. I am also very grateful to Th\'eo Dumont, for our discussions around matrix-valued optimal transport.

\noindent \emph{Use of AI tools.} Generative AI tools (Claude, Anthropic, and ChatGPT, OpenAI) were used during the preparation of this article for literature search, for suggesting and checking some mathematical arguments and computations, for formatting the bibliography, and for editing the manuscript. The author has checked all mathematical content and the accuracy of all references. The author assumes responsibility for all content.

\medskip

\noindent \textbf{Conflict of interest} The author declares that there is no conflict of interest.

\bibliographystyle{siamplain}
\bibliography{biblio}

\end{document}

%% file: introduction.tex
\section*{Introduction}

Convex roofs provide a natural way of extending a function defined on the extreme points of a convex set to the whole set~\cite{uhlmann2010roofs,yan2012extension,bucicovschi2010continuity}. They are constructed with tools from Choquet theory, which connects convex sets with probability measures on their extreme points~\cite{phelps2002lectures,alfsen2012compact}. They have well-established properties, thanks notably to the works of Shirokov~\cite{shirokov2007strong,shirokov2012stability}, see also~\cite{papadopoulou1977geometry,o1976openness}.

\medskip

\noindent Convex roofs are widely used in quantum information theory: quantities defined on pure states, such as the entanglement of a bipartite pure state, are extended to mixed states through their convex roof, yielding entanglement measures~\cite{vollbrecht2001entanglement,uhlmann2010roofs} and, more generally, resource measures~\cite{chitambar2019quantum,zhu2025unified}. These quantities are scalar and, in the developing context of quantum optimal transport~\cite{golse2016mean,golse2017schrodinger,carlen2014analog,de2021quantum,cole2023quantum,toth2023quantum,beatty2025order,friedland2022quantum,agredo2013wasserstein,trevisan2025quantum,beatty2025wasserstein}, it is natural to extend these concepts to functions of two variables, and in particular to distances.

\medskip

\noindent In this work, we therefore study the metric aspects of convex extensions and their relationships with optimal transport on classical and quantum state spaces. We consider convex extensions of functions defined on the extreme points of products of convex sets $C_1 \times C_2$, and specialize to the case where the function to extend is a distance $d$ on the extreme points of a convex set $C$. We construct our extensions as $p$-th roots of the convex extension of $d^p$, $p \in [1,\infty)$. As the triangle inequality may not survive this extension process, we then restore it by taking the metric envelope, obtaining a (pseudo-)distance $D_p$ on $C$. We establish results on the metric structure that it endows $C$ with, gathered in Theorem~\ref{theo:metricDp}.

\noindent We then establish a two-way bridge between convex extensions and optimal transport. The starting point is Choquet theory, which lifts each point of a convex set $C$ to the probability measures on its extreme points ${\rm Ext}\,C$ having this point as barycenter. Conversely, the resultant, or barycenter, map $r : \Pc({\rm Ext}\,C) \to C$ sends each such probability measure back to its barycenter, so that, under suitable assumptions on $C$, the points of $C$ correspond exactly to the classes of probability measures on ${\rm Ext}\,C$ sharing the same barycenter. In other words, any such convex set is in relation, through its resultant map, with a set of probability measures.

\smallskip

\noindent On the one hand, this allows us to recover our convex extensions as quotients of standard, classical Kantorovich costs: the convex extension of a cost between two points of $C$ is obtained by solving the classical Kantorovich problem between their representing measures, and minimizing over all such representations (see Propositions~\ref{prop:foldedKanto} and~\ref{prop:foldedWass}). For this reason, and for brevity, we also refer to our constructions as \emph{folded Kantorovich} and \emph{folded Wasserstein}. In turn, results from optimal transport transfer to general convex sets: for instance, the geodesic structure of Wasserstein spaces transfers to our case.

\smallskip

\noindent On the other hand, various instances of Kantorovich optimal transport formulations re-write as such convex extensions. Geometrically, sets of probability measures and simplices are two sides of the same coin: $\Pc(E)$ can be seen as a simplex whose extreme points, the Dirac masses, correspond to the points of $E$, and conversely, by Bauer~\cite{bauer1961silovscher}, a compact metrizable simplex with closed set of extreme points is affinely homeomorphic to the set of Radon probability measures on its extreme points. This is the case where every point has a unique representation as a barycenter of extreme points (the resultant map is a bijection). Classical optimal transport is in fact the geometrically simplest case of convex extensions, that of the simplex, as already observed by Savar\'e and Sodini~\cite{savare2022simple}. Quantum state spaces, sets of density matrices, are also convex sets, whose extreme points are the pure states, but they are not simplices. Various formulations of quantum optimal transport are framed as problems of extending costs or distances defined between pure states to mixed states, that is from the extreme points of the convex set to the whole set, which is precisely the aim of our construction. We indeed recover the entanglement-free quantum Wasserstein objects of Beatty and Stilck Fran\c{c}a~\cite{beatty2025order} (see also T\'oth and Pitrik~\cite{toth2023quantum,toth2025quantum}). In this non-simplex case, a mixed state admits many representations as a barycenter of pure states, and this degeneracy of representations may break the triangle inequality, hence the need for considering the metric envelope. In this setting, $D_p$ is an actual distance on quantum states as soon as $d$ dominates a norm, for instance for the Frobenius and Fubini--Study distances between pure states (Theorem~\ref{theo:mainquantum} and Corollary~\ref{cor:quantumExamples}). For the Frobenius distance and $p=2$, Miller~\cite{miller2026metricity} recently proved, through a convex-roof duality argument, that the semi-distance of Beatty and Stilck Fran\c{c}a already satisfies the triangle inequality, so that it coincides with $D_2$. The semiclassical optimal transport cost of Golse and Paul~\cite{golse2017schrodinger,golse2021semiclassical,golse2022quantitative}, between quantum and classical states, is also recovered as a convex extension on the product of a quantum and a classical state space. Our framework allows one to transfer tools from the theory of convex extensions to these Kantorovich formulations: notably, the strong CE-property of Shirokov~\cite{shirokov2007strong} yields their continuity.

\smallskip

\noindent In the quantum setting, our framework does not account for entanglement, as it corresponds to optimal transport restricted to separable couplings, as in~\cite{beatty2025order,toth2023quantum,toth2025quantum}. It therefore differs from entanglement-aware formulations, in which the optimization runs over all couplings between quantum states, such as the ones of Golse, Mouhot and Paul~\cite{golse2016mean}, of De~Palma and Trevisan~\cite{de2021quantum} and of Cole, Eckstein, Friedland and {\.Z}yczkowski~\cite{cole2023quantum,friedland2022quantum}. Combining our framework with the noncommutative Choquet theory of Davidson and Kennedy~\cite{davidson2015choquet,davidson2019noncommutative,davidson2024survey}, see also~\cite{kennedy2022noncommutative}, in order to take entanglement into account, is a possible future research direction. Our framework also differs from dynamical and gradient-flow approaches~\cite{carlen2014analog,carlen2017gradient,carlen2020noncommutative} and from matrix-valued optimal transport~\cite{brenier2020matrix}. We refer to the surveys~\cite{trevisan2025quantum,beatty2025wasserstein} for a broader overview. A construction closer in spirit is that of Antonini and Cavalletti~\cite{antonini2021geometry}, whose transport cost between quantum states only involves orthogonal (spectral) decompositions into pure states, instead of all of them as in our construction. This cost is a semi-distance, while fixing the spectral decomposition without repeated eigenvalues removes the degeneracy and yields a distance, possibly infinite. Our framework also differs from weak and barycentric optimal transport~\cite{gozlan2017kantorovich}, as well as from the Choquet theory of Lipschitz-free spaces~\cite{aliaga2024solution}, in which a metric is extended to a norm on a larger linear space. Wasserstein-type distances on Gaussian mixtures~\cite{delon2020wasserstein} and on generic mixture models~\cite{dusson2026wasserstein} can be seen as simplex instances of our construction, with mixtures as convex set and Gaussians (or atoms) as extreme points, although the latter do not (necessarily) form a closed set.

\medskip

\noindent Finally, our framework applies to other convex sets that we do not explore here. The extreme points of the Birkhoff polytope of doubly stochastic matrices are the permutation matrices, so that natural distances between permutations, such as the Cayley or Kendall-$\tau$ distances, can be extended into distances between doubly stochastic matrices. Transportation polytopes of couplings between two discrete measures can be treated similarly. State spaces of general probabilistic theories~\cite{barrett2007information,plavala2023general} are another example. These theories describe a physical system by a convex set of states whose extreme points are the pure states, classical and quantum theories corresponding to simplices and to sets of density matrices. Extending a distance between pure states to all states is then exactly the setting of our construction. The convex set of quantum channels~\cite{choi1975completely} is also a natural candidate, although it is not directly covered, as its extreme points do not form a closed set.

\medskip

\noindent\textbf{Structure of the article.} Section~\ref{sec:ChoquetWasserstein} is devoted to the definition of the core objects of this article, namely the $p$-convex roof extension $\pRoof{p}$ and the $p$-metric convex roof sub-extension $D_p$ of a distance $d$, along with their links with optimal transport and our main Theorem~\ref{theo:metricDp}. In Section~\ref{sec:proofs}, we develop the study of $\pRoof{p}$ and $D_p$ and prove the results stated in Section~\ref{sec:ChoquetWasserstein}. Finally, Section~\ref{sec:applications} is devoted to the applications to classical, semiclassical and quantum structures, and is self-contained assuming the results of Section~\ref{sec:ChoquetWasserstein}.

\medskip

\noindent For a topological space $E$, we denote by $\Pc(E)$ the set of Radon probability measures over $E$. For a topological vector space $X$, we denote by $X^*$ the set of continuous real linear forms over $X$.

\medskip

\noindent \textbf{Pseudo/semi-distance terminology.} A \emph{pseudo}-distance would be distance if separated points. A \emph{semi}-distance would be distance if it satisfied the triangle inequality.

%% file: Section_convex_roofs.tex
\section{Convex roofs and metric envelopes} \label{sec:ChoquetWasserstein} 

This section is devoted to the definition of the core concepts of this article, along with the main results. We start by recalling the notion of \emph{convex roof}, before introducing the \emph{$p$-convex roof extension} and the \emph{$p$-metric convex roof sub-extension} of a distance, of which we then describe the induced metric properties.

\medskip

\noindent Let $K$ be a convex subset of a locally convex Hausdorff vector space $\X$. Its set of extreme points is
\begin{equation}
\label{eqdef:extremepoints}
{\rm Ext}K := \left\{x \in K \; : \; \forall (y,z,t) \in K \times K \times (0,1), \ x = (1-t) y + t z \implies x=y=z \right\}.
\end{equation}
The idea of representing points in $K$ as convex combinations of elements of ${\rm Ext}K$ is formalized, in the framework of Choquet theory, with the following definition.
\begin{definition}[\cite{alfsen2012compact,peifer2020choquet,phelps2002lectures}]
\label{def:representants} A Radon probability measure $\pi \in \Pc({\rm Ext}K)$ is said to \emph{represent} $x \in K$ if for every $f$ in $\X^*$,
\begin{equation}
\label{eq:defRepresentant}
\int_{{\rm Ext}K} f(y) \, \dd \pi(y) = f(x).
\end{equation}

\noindent Conversely, $x$ is called \emph{resultant}, or \emph{barycenter} of any $\pi$ representing $x$. When it can be defined, the application
\begin{equation} \label{eqdef:resultantmap}
    r : \Pc({\rm Ext}K) \to K,
\end{equation}
mapping any $\pi \in \Pc({\rm Ext}K)$ to its barycenter in $K$, is called the \emph{resultant map}. Then $r^{-1}(\{x\})$ is the set of representing measures of $x \in K$ and $r(\pi)$ is the barycenter of $\pi \in \Pc({\rm Ext} K)$.
\end{definition}

\noindent The condition~\eqref{eq:defRepresentant} should simply be understood as the barycenter property
\[
\int_{{\rm Ext} K} y\, \dd \pi(y) = x
\]
without bothering to define an integral on ${\rm Ext} K$.

\medskip

\noindent A sufficient condition for the resultant map to be well-defined is to further assume $K$ and ${\rm Ext}K$ compact~\cite[Proposition~2.5]{peifer2020choquet}, in which case $r$ is even affine and continuous with respect to the weak-$^*$ topology on $\Pc({\rm Ext}K)$~\cite[Lemma~2.6]{peifer2020choquet}. Furthermore, when $K$ is compact metrizable, Choquet's theorem~\cite[Theorem 3.9]{peifer2020choquet} ensures the existence of representing measures for all points in $K$, making $r$ a surjection.

\medskip

\noindent Finally, $K$ is said to be \emph{stable} when the midpoint map $(x,y) \in K \times K \mapsto (x+y)/2 \in K$ is open~\cite{papadopoulou1977geometry,shirokov2012stability}. In the compact setting, stability is equivalent to the openness of the resultant map $r:\Pc(\Ext K)\to K$~\cite{o1976openness}. This notion is fundamental to later obtain continuity properties, notably thanks to the works of Shirokov~\cite{shirokov2012stability}. We highlight that in the compact setting, classical~\cite[Theorem~1]{o1976openness} and quantum~\cite[Theorem~2]{shirokov2006properties} state spaces are indeed stable.

\begin{definition}[\cite{shirokov2007strong}]
\label{def:maximalExtension}Let $K$ be a convex subset of a locally convex Hausdorff vector space $\X$, whose associated resultant map $r$ is well-defined and surjective. Given a lower-bounded measurable $g : {\rm Ext}K \to \R \cup \{+\infty\}$, we define its \emph{maximal convex extension}, or \emph{convex roof}, to be
\begin{equation} \label{eqdef:maxconvexExtension}
\forall x \in K, \qquad     g_*(x) := \inf_{\pi \in r^{-1}(\{x\})} \int_{{\rm Ext}K} g(y) \, \dd \pi(y).
\end{equation}
\end{definition}
\noindent When $K$ and ${\rm Ext} K$ are compact and $g$ is real-valued, lower-bounded and lower semicontinuous (lsc), $g_*$ is the biggest lsc convex function that is upper-bounded by $x \in K \mapsto \begin{cases}
    g(x) \text{ if }  x \in {\rm Ext}K, \\
    + \infty \text{ else,}
\end{cases}$ (hence the naming ``roof''), see~\cite[Theorem~2\,(A)]{shirokov2007strong}.

\medskip

\noindent We now recall the definition of the classical Kantorovich optimal transport cost, later related to convex roofs.

\begin{definition}[\cite{villani2008optimal}]\label{def:optimalTransport}
Consider two Polish spaces $E_1$ and $E_2$ and a lower-bounded measurable \(c : E_1 \times E_2 \to \R \cup \{+\infty\}\).
The \emph{standard Kantorovich} cost associated with $c$ is 
\begin{equation} \label{eqdef:Kantocost}
\forall (\mu,\nu) \in \Pc(E_1) \times \Pc(E_2), \qquad  \mathcal{K}_c (\mu,\nu) = \inf_{\pi \in \Pc_{(\mu,\nu)}} \; \iint_{E_1 \times E_2} c(x,y) \, \dd \pi(x,y),
\end{equation}
with $\Pc_{(\mu,\nu)}$ the set of Radon probability measures over $E_1 \times E_2$ whose marginals are respectively $\mu$ and $\nu$.
\end{definition}

\noindent Specializing Definition~\ref{def:maximalExtension} to the case $K \equiv C_1 \times C_2$ and $g \equiv c : {\rm Ext} C_1 \times {\rm Ext} C_2 \to \R\cup\{+\infty\}$, we can link convex roofs with quotients of Kantorovich costs, relying on the identification
\[
{\rm Ext}(C_1 \times C_2) \;=\; {\rm Ext}\,C_1 \times {\rm Ext}\,C_2.
\]
In turn, such $c_*$ may be called ``folded Kantorovich'' costs.

\begin{proposition}[Folded Kantorovich] \label{prop:foldedKanto}
Let $C_1$ and $C_2$ be two compact metrizable and stable convex subsets of $X_1$ and $X_2$ locally convex Hausdorff vector spaces and $c : {\rm Ext} C_1 \times {\rm Ext} C_2 \to \R\cup\{+\infty\}$ lower semicontinuous. Let $\K_c : \Pc(\Ext C_1) \times \Pc(\Ext C_2) \to \R \cup \{+\infty\}$ be the Kantorovich cost associated with $c$ and denote the resultant maps $r_1 : \Pc(\Ext C_1) \to C_1$ and $r_2 : \Pc(\Ext C_2) \to C_2$. Then the convex roof $c_* : C_1 \times C_2 \to \R \cup \{+\infty\}$ can be recast as the quotient
\begin{equation} \label{eq:foldedKantoRelation}
 c_* = \faktor{\mathcal{K}_c}{(r_1,r_2)},
\end{equation}
in the sense that for any $(x_1,x_2) \in C_1 \times C_2$,
\begin{equation} \label{eq:linkConvexExtensionKantorovich}
c_*(x_1,x_2) = \inf_{\mu_1 \in r_1^{-1}(\{x_1\})} \, \inf_{\mu_2 \in r_2^{-1}(\{x_2\})} \mathcal{K}_c(\mu_1,\mu_2).
\end{equation}
\end{proposition}

\noindent The proof of the above proposition is straightforward from the fact that for any $(x_1,x_2) \in C_1 \times C_2$,
\[
r_{C_1 \times C_2}^{-1}(\{(x_1,x_2)\}) = \bigcup_{(\mu_1,\mu_2) \in r_1^{-1}(\{x_1\}) \times  r_2^{-1}(\{x_2\})} \Pc_{(\mu_1,\mu_2)}.
\]

\noindent We note that the assumptions of Proposition~\ref{prop:foldedKanto} ensure that the ${\rm Ext}\,C_i$ are closed, hence Polish, and the resultant maps $r_1 : \Pc({\rm Ext} C_1) \to C_1$ and $r_2 : \Pc({\rm Ext} C_2) \to C_2$ are well-defined and surjective.

\medskip

\noindent Now specializing further to the case $K = C \times C$ with a single convex set $C$, and the choice $c = d^p$ with $p \in [1,\infty)$ and $d$ a \emph{distance} on ${\rm Ext }C$, we propose the following definition, which is the main object of this work.

\begin{definition}[$p$-metric convex roof sub-extension]
\label{def:metricRoof}Let $C$ be a convex subset of a locally convex Hausdorff vector space $X$, whose associated resultant map is well-defined and surjective. Let $p \in [1,\infty)$ and $d : {\rm Ext}C \times {\rm Ext}C \to \R_+$ a measurable distance on ${\rm Ext}C$. We define the \emph{$p$-convex roof extension of $d$} as the application
\begin{equation} \label{eqdef:PhiMaxConvExt}
    \pRoof{p} := ((d^p)_*)^{\frac1p} : C \times C \to [0,+\infty].
\end{equation}
The \emph{$p$-metric convex roof sub-extension} of $d$ is then defined to be the metric envelope of $\pRoof{p}$, that is, for any $(x,y) \in C \times  C$,
\begin{equation} \label{eqdef:PhimaxMetricconvexExtension}
D_{p}(x,y) := \inf\left\{\sum_{j=1}^N \pRoof{p}(z_j,z_{j+1}) \; : \; N \in \N, \ z_1 = x, \ z_{N+1} = y, \ z_j \in C \right\}.
\end{equation}
\end{definition}

\noindent Note that, despite the terminology, $(\pRoof{p})^p$ is indeed convex as a convex roof, but $\pRoof{p}$ is not in general for $p>1$. A priori, $\pRoof p$ may not preserve the triangle inequality, and the idea of~\eqref{eqdef:PhimaxMetricconvexExtension}, standard when defining pseudo-distances on quotient spaces, is to enforce the subadditivity by taking the minimal length among all possible chains. This procedure however does not a priori preserve the separation property. $D_p$, the metric envelope of $\pRoof p$, is the largest pseudo-distance that is no more than $\pRoof p$, see for instance~\cite{burago2001course}. In particular, if $\pRoof p$ is a distance, then $\pRoof p = D_p$. In general, $\pRoof{p}$, and hence $D_p$, may take the value $+\infty$, but they are finite in the case where $C$ is compact and $d$ continuous. The terminology is justified by Proposition~\ref{prop:extensionproperty}: $\pRoof{p}$ extends $d$, while $D_p$ in general only sub-extends it, in the sense that ${D_p}_{|\Ext C\times\Ext C}\leq d$. For short, and in view of Proposition~\ref{prop:foldedKanto} and the next Proposition~\ref{prop:foldedWass}, $\pRoof{p}$ and $D_p$ are respectively referred to as \emph{folded Kantorovich-$p$} and \emph{folded Wasserstein-$p$}. We recall that the \emph{standard Wasserstein}-$p$ distance associated with a Polish space $(E,d)$ is defined from the Kantorovich cost on $\Pc(E)\times\Pc(E)$ by
\begin{equation} \label{eqdef:Wp}
    W_p := \left(\mathcal{K}_{d^p}\right)^{\frac1p}.
\end{equation}

\begin{proposition}[Folded Wasserstein] \label{prop:foldedWass}
Let $C$ be a \emph{compact} and \emph{stable} convex subset of a locally convex Hausdorff vector space and let $d$ be a distance on ${\rm Ext}\,C$ continuous on $\Ext C\times\Ext C$ with respect to the natural topology. Denote by $r$ the resultant map for the convex set $C$. Then for any $p \in [1,\infty)$, $D_p$ is the \emph{quotient pseudo-distance} of the Wasserstein-$p$ distance $W_p$ on $\Pc({\rm Ext} C)$ associated with $d$, by the surjection $r : \Pc({\rm Ext} C) \to C$, that is~\cite[Definition 3.1.12]{burago2001course}\footnote{ In this reference, the terminology ``semi-metric'' is used for what is denoted here ``pseudo-distance''.} the largest pseudo-distance on $C$ that is no more than
\begin{equation*}
\pRoof{p} : (x,y) \in C \times C \mapsto \inf_{\mu \in r^{-1}(\{x\})} \, \inf_{\nu \in r^{-1}(\{y\})} W_p(\mu,\nu).
\end{equation*}
\end{proposition}

\noindent Note that, $\Ext C$ being compact under these assumptions, such a continuous distance $d$ automatically induces the natural topology on $\Ext C$ (see Section~\ref{sec:proofs}).

\medskip

\noindent We finally state our main Theorem~\ref{theo:metricDp}, exposing the metric properties $D_p$ gives on $C$, under an assumption~\eqref{eq:sepHyp} ensuring it to separate the points of $C$. This assumption is satisfied in particular when $X$ is a Banach space whose norm distance $\|-\cdot -\|$ is dominated by $d$ on $\Ext C \times \Ext C$.
\begin{theorem}\label{theo:metricDp}
Let $C$ be a \emph{compact} and \emph{stable} convex subset of a locally convex real Hausdorff vector space~$X$. Let $d$ be a distance on ${\rm Ext} C$ continuous on $\Ext C \times \Ext C$ with respect to the natural topology. Further assume that
\begin{equation} \label{eq:sepHyp}
\text{the real continuous affine } d\text{-Lipschitz functions over } \Ext C \text{ separate the points of } C.
\end{equation}
Then $(D_p)_{p\geq 1}$ is a nondecreasing family of distances on $C$, whose induced topologies all coincide with the natural one, and in particular the $(C,D_p)$ are compact. 

\smallskip

\noindent Moreover, $(C,D_1)$ is always geodesic, and if $({\rm Ext} C, d)$ is geodesic, then $(C, \, D_p)$ is also geodesic for all $p>1$.

\smallskip

\noindent Finally, for all $p \in [1,\infty)$,  $D_p$ \emph{sub}-extends $d$, in the sense that 
$$
{D_p}_{|{\rm Ext}C\times{\rm Ext}C} \leq d.$$
If moreover $X \equiv (X,\|\cdot\|)$ is a Banach space endowed with its norm topology and $d(z,z_*)=\|z-z_*\|$ for all $z,z_*\in\Ext C$, then $D_p$ extends $d$, in the sense that 
\begin{equation} \label{eqTheo:ExtendsYes}
    {D_p}_{|{\rm Ext}C\times{\rm Ext}C} = d. 
\end{equation}
\end{theorem}

\noindent The proof of Theorem~\ref{theo:metricDp} can be found at the very end of Section~\ref{sec:proofs} after all our developments around $D_p$. While Assumption~\eqref{eq:sepHyp} may appear difficult to test, it holds in various cases of interest for us.

\begin{itemize}
\item[\emph{(i)}] (Classical state space) If $C = \Pc(E_0)$ for some compact metric space $(E_0,d_0)$. Then $\Ext C = \{\delta_x\}_{x \in E_0}$ can be naturally endowed with $d(\delta_x,\delta_y) = d_0(x,y)$. Then the set of continuous affine $d$-Lipschitz functions over $\Ext C$ is isomorphic to the set of $d_0$-Lipschitz functions over $E_0$. The fact that these functions separate the set of probability measures over $E_0$ can be recovered from the Kantorovich--Rubinstein duality for the Wasserstein-1 distance~\cite{villani2008optimal}. Note finally that in this case, we simply recover
\begin{equation*}
    D_p = W_p,
\end{equation*}
with $W_p$ the Wasserstein-$p$ distance defined from $(E_0,d_0)$.
\item[\emph{(ii)}] (Norm domination) If $X\equiv(X,\|\cdot\|)$ is a Banach space and $d$ upper-bounds the norm, in the sense that $d(x,y) \geq \|x-y\|$ for all $x,y \in \Ext C$, then the set of continuous linear forms over $X$ with dual norm one is contained in the set of continuous affine $d$-Lipschitz functions over $\Ext C$. Then~\eqref{eq:sepHyp} comes by duality of the norm, for any $x \in X$, $\|x\| = \underset{f \in X^*, \, \|f\|_{X^*} = 1}{\sup} f(x)$.

\smallskip

\noindent This norm domination condition is satisfied (for some suitable norm) for various relevant distances in the case of quantum state spaces (like the Fubini--Study distance).
\end{itemize}

\begin{remark}[Dual formulation]
A Fenchel--Moreau~\cite{rockafellar1998variational,volle2015duality} argument can provide a dual formulation for~$c_*$, and hence $\pRoof{p}$, analogous to the Kantorovich--Rubinstein duality and generalizing~\cite[Theorem~4.4]{savare2022simple} which derives it in the case where the convex sets are simplexes. We did not investigate further this duality in this work. In the quantum setting with the Hilbert--Schmidt cost and $p=2$, such a dual formulation is derived and exploited in~\cite{miller2026metricity}.
\end{remark}

%% file: Section_proofs.tex
\section{Developments on the \texorpdfstring{$p$}{p}-metric convex roof sub-extension and proof of \texorpdfstring{Theorem~\ref{theo:metricDp}}{the main theorem}} \label{sec:proofs}

This section is devoted to the development and study of the folded Wasserstein-$p$ pseudo-distance $D_p$. Throughout, $p \in [1,\infty)$, $X$ denotes a locally convex real Hausdorff vector space, whose topology is called the \emph{natural topology}, $C$ is a nonempty compact and stable convex subset of $X$ and $\Ext C$ is the set of extreme points of $C$ (its Choquet boundary). We consider $d$ a distance on $\Ext C$ continuous on $\Ext C\times\Ext C$ with respect to the natural topology. Since the extreme boundary of a stable convex set is always closed~\cite{papadopoulou1977geometry}, ${\rm Ext}\,C$ is compact and therefore $d$, being continuous, induces the natural topology on $\Ext C$, and $(\Ext C,d)$ is a compact, hence Polish, metric space. This moreover ensures $C$ to be metrizable~\cite{corson1970metrizability}. All these properties of $C$ (compactness, stability, metrizability, closedness of the extreme boundary) are inherited by $C\times C$, whose extreme boundary is $\Ext C\times\Ext C$.

\medskip

We denote by $r : \Pc(\Ext C) \to C$ the resultant map, which our assumptions ensure to be well-defined, affine and continuous~\cite[Proposition~2.5]{peifer2020choquet}, and, by stability, open~\cite{o1976openness}. Informally, it sends a convex combination of points on $\Ext C$ to its corresponding barycenter in $C$. By the Krein--Milman theorem~\cite[Proposition~1.2]{phelps2002lectures} (or Choquet's theorem~\cite[Section~3]{phelps2002lectures}), our assumptions ensure that every point in $C$ has a representing measure in $\Pc(\Ext C)$. As such, $r$ is surjective, and this allows one to identify $C$ with the corresponding quotient,
\begin{equation} \label{eq:ChoquetRepresentation}
    C \; \; \cong \; \;\faktor{\Pc (\Ext C)}{r},
\end{equation}
where the quotient by $r$ identifies $\mu,\nu \in \Pc(\Ext C)$ whenever they have the same barycenter, $r(\mu)=r(\nu)$.
For any $x,y \in C$, we denote $\Pc_x := r^{-1}(\{x\})$ and $\Pc_{(x,y)} := r_{C \times C}^{-1}\{(x,y)\}$ (with $r_{C \times C}$ the resultant map associated with the convex set $C \times C$). In particular, 
\[
\Pc(\Ext C) = \bigsqcup_{x \in C} \Pc_x, \qquad \Pc(\Ext C \times \Ext C) = \bigsqcup_{(x,y) \in C \times C} \Pc_{(x,y)}.
\]
 Since $\Pc(\Ext C)$ is weak-$*$ compact, $r$ (resp. $r_{C \times C}$) is continuous from the weak-$*$ topology to the natural topology~\cite[Lemma~2.6]{peifer2020choquet} and $\{x\}$ is closed, the fibers $\Pc_x$ and $\Pc_{(x,y)}$ are compact.

\medskip
 
One finding of the present work is (see Proposition~\ref{prop:foldedWass}) that $D_p$, the metric envelope of~$\pRoof{p}$, is equal to the quotient (pseudo-)distance of the standard Wasserstein\nobreakdash-$p$~distance $W_p$ by the barycenter equivalence, in the sense of~\cite[Definition 3.1.12]{burago2001course}, given by the process
\begin{equation} \label{eq:foldingProcess}
(\Ext C,d) \underset{\text{Standard OT}}{\longrightarrow} \left(\Pc(\Ext C),W_p \right) \underset{\text{Gluing}}{\longrightarrow} \left(\faktor{\Pc (\Ext C)}{r}, \faktor{W_p}{r} \right) \cong (C,D_p).
\end{equation}
Indeed, combining Proposition~\ref{prop:foldedKanto} with~\eqref{eqdef:Wp} one finds that for all $x,y \in C$,
\begin{equation} \label{eq:pRoofWp}
    \pRoof{p}(x,y) = \inf_{\mu \in \Pc_x} \, \inf_{\nu \in \Pc_y} W_p(\mu,\nu),
\end{equation}
of which $D_p$ is the metric envelope. We now turn to studying the properties of the folded Kantorovich-$p$ semi-distance $\pRoof{p}$ and the folded Wasserstein-$p$ pseudo-distance $D_p$, which altogether contribute to the proof of our main Theorem~\ref{theo:metricDp}. We start with basic properties.

\begin{proposition} \label{prop:propertiesDp} $D_p$ is a pseudo-distance on $C$ and $\pRoof{p}$ is a semi-distance on $C$, with
\begin{equation} \label{eqprop:comparisonQuotient}
D_p \leq \pRoof{p}.
\end{equation} 
Moreover,
\begin{equation} \label{eqprop:equivQuotient}
    \pRoof{p}\; \; \text{is a distance} \; \; \iff  \pRoof{p}\; \; \text{is subadditive} \; \; \iff  \; \; \pRoof{p} = D_p.
\end{equation}
\end{proposition}
\begin{proof}
This whole proposition (both~\eqref{eqprop:comparisonQuotient} and~\eqref{eqprop:equivQuotient}) is an immediate consequence of $W_p$ being a distance on $\Pc(\Ext C)$ and $D_p$ being the quotient pseudo-distance~\cite[Definition 3.1.12]{burago2001course} on $\faktor{\Pc(\Ext C)}{r} \cong C$. The only nontrivial matter to prove is the separation property of $\pRoof{p}$.

\smallskip

\noindent For any $x \neq y$, $\Pc_x$ and $\Pc_{y}$ (see Definition~\ref{def:representants}) are non-empty, disjoint (as equivalence classes) and closed in the weak-$^*$ natural topology, by continuity of the barycenter map. As noted at the beginning of this section, $d$ induces the natural topology, and by~\cite[Theorem~6.9]{villani2008optimal} $W_p$ metrizes the weak-$^*$ convergence for the topology induced by $d$. Moreover, the compactness of $(\Ext C,d)$ implies the one of $(\Pc(\Ext C), \, W_p)$, see~\cite[Remark~6.19]{villani2008optimal}. 

\smallskip

\noindent Therefore $\Pc_x$ and $\Pc_y$ are non-empty, disjoint and compact sets of $(\Pc(\Ext C), \, W_p)$, which implies using~\eqref{eq:pRoofWp}
\[
\pRoof{p}(x,y) = \inf_{\mu \in \Pc_x} \, \inf_{\nu \in \Pc_{y}} \, W_p(\mu,\nu) > 0,
\]
meaning that $\pRoof{p}$ separates points.
\end{proof}

\noindent Since $\pRoof{p}$ is easier to manipulate than $D_p$, we provide a sufficient condition for it to be a genuine distance.

\begin{proposition} \label{prop:subadditivitygeneral}
Assume that for every $x,y,z\in C$ there exists $\mu^*\in\Pc_x$ realizing the two infima
\begin{equation} \label{eq:commonrepresentative}
\pRoof{p}(x,y) = \inf_{\nu \in \Pc_y} W_p(\mu^*,\nu) \qquad \text{and} \qquad \pRoof{p}(x,z) = \inf_{\nu \in \Pc_z} W_p(\mu^*,\nu).
\end{equation}
Then $\pRoof{p}$ is subadditive.
\end{proposition}
\begin{proof}
Take some triplet $x,y,z \in C$ and the associated $\mu^*$. Then for any $(\nu, \eta) \in \Pc_y \times \Pc_z$, we have from~\eqref{eq:pRoofWp} and the subadditivity of $W_p$ that
\[
\pRoof{p}(y,z) \leq  W_p(\nu,\eta) \leq W_p(\nu,\mu^*) + W_p(\mu^*,\eta).
\]
Taking the infimum over $\nu$ and $\eta$ on the right-hand side while using~\eqref{eq:commonrepresentative} yields
\[
\pRoof{p}(y,z) \leq \pRoof{p}(y,x) + \pRoof{p}(x,z).
\]
\end{proof}

\begin{proposition}\label{prop:extensionproperty}
 On $C$, $\pRoof{p}$ extends $d$ while $D_p$ sub-extends $d$, in the sense that
\begin{equation} \label{eq:DpExtendsd}
{{}\pRoof{p}}_{|\Ext C\times \Ext C} = d \qquad \text{and} \qquad {D_p}_{|\Ext C\times \Ext C} \leq d.
\end{equation}
\end{proposition}
\begin{proof}
By~\cite[Proposition~1.4]{phelps2002lectures} (Bauer~\cite{bauer1961silovscher}), extreme points are represented solely by their associated Dirac mass, i.e. for any $x \in \Ext C$, the set $\Pc_x$ of representing probabilities of $x$ is  exactly $\{\delta_x\}$. As such, and using~\eqref{eq:pRoofWp}, we have for all $(x,y) \in \Ext C \times \Ext C$ that
\[
\pRoof{p}(x,y) = \inf_{\mu \in \Pc_x} \, \inf_{\nu \in \Pc_y} \, W_p(\mu, \nu) = W_p(\delta_x, \delta_y) = d(x,y),
\]
and we conclude with~\eqref{eqprop:comparisonQuotient}.
\end{proof}

\smallskip

\noindent The next proposition is a generalization of the same monotonicity property for the standard Wasserstein family of distances~\cite[Remark~6.6]{villani2008optimal}.

\begin{proposition}\label{prop:monotonicity}
 For every $q \in [p,\infty)$, we have the monotonicity properties
\begin{equation} \label{eqprop:DplessDq}
    \pRoof{p} \leq \pRoof{q} \leq \mathrm{diam}_d(\Ext C)^{1 - \frac{p}{q}} \; (\pRoof{p})^\frac{p}{q} \qandq D_p \leq D_q,
\end{equation}
with $\mathrm{diam}_d(\Ext C) = \sup_{x,y \in \Ext C} d(x,y)$ the diameter of $\Ext C$.
\end{proposition}
\begin{proof}
Let $x,y \in C$ and a plan $\pi \in \Pc_{(x,y)}$. By Jensen's inequality we have
\[
\left(\iint_{\Ext C \times \Ext C } d(z,z_*)^p\,\dd\pi(z,z_*) \right)^{1/p} \leq \left(\iint_{\Ext C \times \Ext C } d(z,z_*)^q\,\dd\pi(z,z_*) \right)^{1/q}.
\]
By taking the infimum over all $\pi  \in \Pc_{(x,y)}$  in the left-hand-side, and then on the right-hand-side of the above inequality, we conclude by definition~\eqref{eqdef:maxconvexExtension} that
\[
\pRoof{p}(x,y) \leq \pRoof{q}(x,y).
\]
Conversely, it is straightforward to note that
\begin{multline*}
\left(\iint_{\Ext C \times \Ext C } d(z,z_*)^q\,\dd\pi(z,z_*) \right)^{1/q} \\
\leq \mathrm{diam} (\Ext C)^{1-p/q}\left(\left(\iint_{\Ext C \times \Ext C } d(z,z_*)^p\,\dd\pi(z,z_*) \right)^{1/p} \right)^{p/q},
\end{multline*}
leading to $\pRoof{q} \leq \mathrm{diam}_d(\Ext C)^{1 - \frac{p}{q}} (\pRoof{p})^{\frac{p}{q}}$.

\smallskip

\noindent Since we just proved that every chain of points of $C$ yields a smaller $\pRoof{p}$-length that $\pRoof{q}$-length, it comes naturally by~\eqref{eqdef:PhimaxMetricconvexExtension} that $D_p \leq D_q$. The converse inequality however may not survive the chaining step.
\end{proof}

The established monotonicity of the families $(\pRoof{q})_{q \geq 1}$ and $(D_q)_{q\geq 1}$ allow us to consider the limiting objects $q \to \infty$. Although not used in the rest of the paper, we give the following

\begin{definition}
The folded Kantorovich\nobreakdash-$\infty$~semi-distance and folded Wasserstein\nobreakdash-$\infty$~pseudo-distance are defined by
    \begin{equation}
        \pRoof{\infty} := \lim_{q \to \infty} \pRoof{q} = \sup_{q \geq 1} \; \pRoof{q}.
    \end{equation}
    and
        \begin{equation}
        D_\infty := \lim_{q \to \infty} D_q = \sup_{q \geq 1} \;  D_q.
    \end{equation}
\end{definition}

\begin{proposition} \label{prop:orderinfty}
On $C$, $\pRoof{\infty}$ is a semi-distance while $D_\infty$ is pseudo-distance. Moreover, letting $\widehat{\pRoof{\infty}}$ be the metric envelope of $\pRoof{\infty}$, we have
\begin{equation} \label{eq:ineqProfinfty}
    D_\infty \leq \widehat{\pRoof{\infty}}.
\end{equation}
\end{proposition}
\begin{proof}
First notice that $\pRoof{\infty}$ is finite since Proposition~\ref{prop:monotonicity} implies that
\[
\pRoof{\infty} \leq \mathrm{diam}_d(\Ext C) < \infty.
\]
This also holds for $D_\infty$ as
\begin{equation*} 
    \forall q\geq 1, \quad D_q \leq \pRoof{q} \leq \pRoof{\infty}
\end{equation*}
implies by taking the supremum that
\begin{equation} \label{eq:proofDinfty}
D_\infty \leq \pRoof{\infty}.
\end{equation}
As a bounded supremum of semi-distances, $\pRoof{\infty}$ is itself a semi-distance, and being a bounded supremum of pseudo-distances, $D_\infty$ is a pseudo-distance. By definition, $\widehat{\pRoof{\infty}}$ is the largest pseudo-distance that does not exceed $\pRoof{\infty}$, and therefore~\eqref{eq:proofDinfty} implies~\eqref{eq:ineqProfinfty}.
\end{proof}

\smallskip

In the next proposition, we show the existence of a minimizer, and provide the existence of a finitely-supported minimizer with an upper-bounded number of atoms in the case where representing probabilities are characterized by a finite number of linear constraints.

\begin{proposition}
\label{prop:linearProgramhatDp} 
For all $(x,y) \in C \times C$, there exists a $\pi^* \in \Pc_{(x,y)}$ realizing
\begin{equation}
\iint_{\Ext C \times \Ext C} \, d(z,z_*)^p \, \dd \pi^*(z,z_*) = \inf_{\pi \in \Pc_{(x,y)}} \iint_{\Ext C \times \Ext C} \, d(z,z_*)^p \, \dd \pi(z,z_*),
\end{equation}
called  \emph{optimal representing transport plan}. If moreover there exist $N_x,N_y \in \N^*$, two families of measurable real functions $\{f^x_i\}_{1\leq i \leq N_x-1}$, $\{f^y_i\}_{1\leq i \leq N_y-1}$ over $\Ext C$ and real scalars $\{c^x_i\}_{1\leq i \leq N_x-1}$, $\{c^y_i\}_{1\leq i \leq N_y-1}$ such that, for $z \in \{x,y\}$,
\begin{equation}
  \Pc_z = \left\{\mu \in \Pc(\Ext C) \; : \; f^z_i \in L^1(\mu) \text{ and } \int_{\Ext C} f^z_i \, \dd \mu = c^z_i \text{ for } 1 \leq i \leq N_z-1 \right\},
\end{equation}
then there exists a finitely supported optimal representing transport plan with at most $N_x+ N_y-1$ atoms.
\end{proposition}
\begin{proof}
Since $\Ext C \times \Ext C$ is compact, $d$ is continuous and bounded, so that the application
\[
G : \pi \in \Pc(\Ext C \times \Ext C) \quad \mapsto \quad \iint_{\Ext C \times \Ext C} \, d(z,z_*)^p \, \dd \pi(z,z_*)
\]
is continuous for the weak-$^*$ topology. As the resultant map $r_{C \times C}$ is continuous, the set $\Pc_{(x,y)}$ is weak-$^*$ closed in the compact $\Pc(\Ext C \times \Ext C)$, hence compact. Therefore the infimum of $G$ over $\Pc_{(x,y)}$ is attained.

\smallskip

\noindent As for the second part of the proposition, by~\cite[Proposition~3.1]{winkler1988extreme}, the application $G$ is measure affine, so that we may apply~\cite[Theorem~3.2 with $H = \Pc_{(x,y)}$ and $F=-G$]{winkler1988extreme}, which ensures that
\[
\inf_{\pi \in \Pc_{(x,y)}} G(\pi) = \inf_{\pi \in \partial_e \Pc_{(x,y)}} G(\pi),
\]
where $\partial_e \Pc_{(x,y)}$ denotes the set of extreme points of the compact convex set $\Pc_{(x,y)}$. Next we notice that
\begin{multline}
\Pc_{(x,y)} = \Big\{\pi \in \Pc(\Ext C \times \Ext C) \; : \; f^{x,y}_k \in L^1(\pi) \text{ and } \\
\iint_{\Ext C \times \Ext C} f^{x,y}_k \, \dd \pi = c^{x,y}_k \text{ for } 1 \leq k \leq N_x + N_y-2 \Big\},
\end{multline}
with, for all $1 \leq k \leq N_x-1$,
\[
\forall z,z_* \in \Ext C, \quad f^{x,y}_k(z,z_*) = f^x_k(z) \quad \text{and} \quad c^{x,y}_k = c^x_k, 
\]
and  for all $N_x \leq k \leq N_x + N_y - 2$,
\[
\forall z,z_* \in \Ext C, \quad  f^{x,y}_{k}(z,z_*) = f^y_{k - N_x + 1}(z_*) \qandq c^{x,y}_{k} = c^y_{k - N_x + 1}.
\]
We may then use~\cite[Theorem~3.1~(b)]{winkler1988extreme} to obtain that
\[
\partial_e \Pc_{(x,y)} \subset \Pc_{\leq N_x + N_y-1}(\Ext C \times \Ext C) \cap \Pc_{(x,y)},
\]
where $\Pc_{\leq N_x + N_y-1}(\Ext C \times \Ext C)$ stands for the Radon probability measures on $\Ext C \times \Ext C$ which have finite support and at most $N_x + N_y-1$ atoms. All in all, we have
\begin{equation}
    \pRoof{p}(x,y)^p = \inf_{\pi \in \Pc_{(x,y)}} G(\pi) = \inf_{\pi \in \Pc_{\leq N_x + N_y-1}(\Ext C \times \Ext C) \cap \Pc_{(x,y)}} G(\pi).
\end{equation}
Since the set $\Pc_{\leq N_x + N_y-1}(\Ext C \times \Ext C)\cap \Pc_{(x,y)}$ is closed in the compact set $\mathcal{P}(\Ext C \times \Ext C)$ and $G$ is continuous, the infimum in the right-hand-side is attained. We conclude to the existence of a finitely supported $\pi^* \in \Pc_{(x,y)}$ with at most $N_x + N_y-1$ atoms such that
\[
\pRoof{p}(x,y)^p = \iint_{\Ext C \times \Ext C} \, d(z,z_*)^p \, \dd \pi^*(z,z_*).
\]
\end{proof}

\noindent In the next proposition, we show that $\pRoof{p}$ and $D_p$ are continuous in the natural topology. This is a direct consequence of $\pRoof{p}^p$ being the convex extension on $C \times C$ of the function $d^p$, continuous on $\Ext C \times \Ext C$, and the fact that the convex set $C$ is assumed to be compact and stable.

\begin{proposition}\label{prop:hatDpContinuousGeneral}
On $C \times C$, the applications $\pRoof{p}$ and $D_p$ are continuous with respect to the natural topology.
\end{proposition}

\begin{proof} Since $C$ is compact and stable, so is $C\times C$, and its extreme boundary $\Ext(C\times C)=\Ext C\times\Ext C$ is compact. The continuous function $d^p$ is therefore bounded on $\Ext C\times\Ext C$. Moreover, $C\times C$ satisfies the strong CE-property, as a stable compact convex set~\cite[$(i) \implies (iii)$ in Theorem~1]{shirokov2007strong}. Since $(\pRoof{p})^p$ is the convex extension of $d^p$ on $C\times C$, we obtain by~\cite[Corollary~3]{shirokov2007strong} that $(\pRoof{p})^p$ is continuous and bounded on $C\times C$, hence so is $\pRoof{p}$.

\medskip

\noindent Now, let us prove the continuity of $D_p$. Let $(x_n,y_n)_n$ be a sequence of elements of $C \times C$ converging in the natural topology towards $(x_\infty,y_\infty) \in C \times C$. Then Proposition~\ref{prop:propertiesDp} and the continuity of $\pRoof{p}$ give
\[
0 \leq D_p(x_n,x_\infty) \leq \pRoof{p}(x_n,x_\infty) \xrightarrow[n\to\infty]{} \pRoof{p}(x_\infty,x_\infty) = 0.
\]
Similarly, $D_p(y_n,y_\infty) \xrightarrow[n\to\infty]{}  0$. Since $D_p$ is subadditive, we have
\[
D_p(x_n,y_n) \leq D_p(x_n,x_\infty) + D_p(x_\infty,y_\infty) + D_p(y_\infty,y_n),
\]
so that $\lim \sup_n D_p(x_n,y_n) \leq D_p(x_\infty,y_\infty)$. Conversely,
\[
D_p(x_\infty,y_\infty) \leq D_p(x_\infty,x_n) + D_p(x_n,y_n) + D_p(y_n,y_\infty),
\]
so that $D_p(x_\infty,y_\infty) \leq \lim \inf_n D_p(x_n,y_n)$. This ends the proof, since \(C \times C\) is a metrizable space, in which sequential continuity implies continuity.
\end{proof}

\noindent We have moreover a convexity property for $\pRoof{1}$ and $D_1$.

\begin{proposition}\label{prop:convexityDp}
In the case $p=1$, both the folded Kantorovich-$1$ semi-distance $\pRoof{1}$ and the folded Wasserstein-$1$ pseudo-distance $D_1$ are jointly convex on $C \times C$.
\end{proposition}
\begin{proof}
By definition~\eqref{def:maximalExtension}--\eqref{eqdef:PhiMaxConvExt}, $\pRoof{1} \equiv d_*$ is the convex roof extension of $d$ onto $C \times C$, and is therefore trivially jointly convex. As of $D_1$, consider $(x_0,x'),(y_0,y') \in C \times C$, $t \in [0,1]$, and let $x_0,\dots,x_{N+1}=x'$ and $y_0,\dots,y_{N+1}=y'$ be two chains in $C$ (that can be taken of the same length, up to repeating points). Then $z_i^t := t x_i + (1-t) y_i$ is a chain from $t x_0 + (1-t) y_0$ to $t x' + (1-t) y'$, so that, using the joint convexity of $\pRoof{1}$, we have
\begin{align*}
D_1\big(t(x_0,x') + (1-t)(y_0,y')\big) &\leq \sum_{i=0}^N \pRoof{1}(z_i^t,z_{i+1}^t) \\
&\leq \sum_{i=0}^N \big[\, t\,\pRoof{1}(x_i,x_{i+1}) + (1-t)\,\pRoof{1}(y_i,y_{i+1}) \,\big].
\end{align*}
Taking the infimum over both $x$- and $y$-chains yields 
$$
D_1\big(t(x_0,x')+(1-t)(y_0,y')\big) \leq t\,D_1(x_0,x') + (1-t)\,D_1(y_0,y').
$$
\end{proof}

\smallskip

\noindent Let us now focus on the structure that $D_p$ gives to the space $C$, when we assume that it separates points. We highlight in particular the ease with which the geodesic property is obtained, and highlight that it encourages the investigation for an equivalent dynamical formulation.

\begin{proposition}\label{prop:structureCDp}
Assume $D_p$ to separate points. Then it is a distance on $C$ whose induced topology is equivalent to the natural topology, and $(C,D_p)$ is compact. Moreover $(C,D_1)$ is always geodesic, and if $(\Ext C, d)$ is geodesic, then $(C,D_p)$ is geodesic for $p>1$.
\end{proposition}
\begin{proof}

Since by Proposition~\ref{prop:propertiesDp}, $D_p$ is a pseudo-distance, it is a distance as soon as it separates points. 

\smallskip

$\bullet$ We start with the compactness and topology part. From~\cite[Theorem~6.18 and Remark~6.19]{villani2008optimal}, the fact that $(\Ext C, \, d)$ is Polish and compact implies that $(\Pc(\Ext C),  \, W_p)$ is also Polish and compact. Since $(C, \, D_p)$ is the quotient space (endowed with the quotient distance) of $(\Pc(\Ext C),  \, W_p)$ for the equivalence relation induced by the resultant map $r$, its topology coincides with the quotient topology~\cite[Exercise~3.1.14]{burago2001course}. Moreover Proposition~\ref{prop:hatDpContinuousGeneral} ensures that $D_p$ is continuous with respect to the natural topology. Therefore the identity map from $C$, endowed with the natural topology, onto $(C,D_p)$ is a continuous bijection from a compact space onto a Hausdorff space, hence a homeomorphism. Therefore the two topologies coincide.

\smallskip

$\bullet$ We now turn to the geodesic part. Let first $p=1$, $x,y\in C$ and $x_t:=(1-t)x+ty$ for $t\in[0,1]$. For $0\leq s\leq t\leq1$, $(x_s,x_t)=(1-t)(x,x)+(t-s)(x,y)+s(y,y)$ is a convex combination in $C\times C$, so the joint convexity of $D_1$ (Proposition~\ref{prop:convexityDp}) yields $D_1(x_s,x_t)\leq(t-s)\,D_1(x,y)$. Hence, by the triangle inequality,
\[
D_1(x,y)\leq D_1(x,x_s)+D_1(x_s,x_t)+D_1(x_t,y)\leq\big(s+(t-s)+(1-t)\big)D_1(x,y)=D_1(x,y),
\]
so that all these inequalities are equalities, and $D_1(x_s,x_t)=(t-s)\,D_1(x,y)$: the segment $t\mapsto x_t$ is a constant-speed geodesic. Now assume $p>1$ and $(\Ext C, \,d)$ geodesic. From~\cite[Corollary~7.22]{villani2008optimal}, $(\Pc(\Ext C), \, W_p)$ is also geodesic. On the other hand, it is ensured in~\cite[Section~3.1]{burago2001course} that the quotient space of a length space is also a length space. Remarking that a compact length space is geodesic, by the Hopf--Rinow--Cohn--Vossen theorem~\cite[Theorem~2.5.28 and Remark~2.5.29]{burago2001course}, we conclude that $(C, \, D_p)$ is geodesic.
\end{proof}

\noindent The next lemma allows one to ensure the separation property of $D_p$ under some assumption on $d$. We mention both that the assumption~\eqref{eq:lowerDbd} corresponds to the one in Beatty~\cite[Equation~(26)]{beatty2025order} and that it is satisfied for the two typical choices of $d$ we make in the Section~\ref{sec:quantumChoquetWasserstein} (see Corollary~\ref{cor:quantumExamples}).

\begin{proposition}\label{prop:sepLipschitz}
Denote by $A_{d\text{-Lip}}(C)$ the real continuous and affine functions on $C$ that are $d$-Lipschitz on $\Ext C$, 
\begin{equation} \label{eqdef:lipSet}
    A_{d\text{-Lip}}(C) := \{a \in A(C), \; \; |a(z)-a(z_*)| \leq d(z,z_*) \; \; \forall z,z_* \in \Ext C\}.
\end{equation}
Then for any $a \in A_{d\text{-Lip}}(C)$ we have
\[ \forall\,x,y\in C,\qquad D_p(x,y) \;\geq\; |a(x)-a(y)|. \]
If in particular $X \equiv (X,\|\cdot\|)$ is a Banach space and
\begin{equation} \label{eq:lowerDbd}
\forall \, z,z_* \in \Ext C, \qquad d(z,z_*) \geq \|z-z_*\|, \phantom{\forall \, z,z_* \in \Ext C, \qquad}
\end{equation}
then
\begin{equation} \label{eq:lowerBdDp}
\forall \, x,y \in C, \qquad D_p(x,y) \geq \|x-y\|. \phantom{\forall \, x,y \in C, \qquad}
\end{equation}
\end{proposition}

\begin{remark}
We can also prove that if $d(z,z_*) \geq \|z-z_*\|^\beta$ with $\beta > 1$, then $\pRoof{p}(x,y) \geq \|x-y\|^\beta$, but nothing can be said of $D_p$ in this case.
\end{remark}

\begin{proof}
The proof of~\eqref{eq:lowerBdDp} is related to the one of~\cite[Lemma~12 with $\alpha = 1$]{beatty2025order}. Since by Proposition~\ref{prop:monotonicity}, $D_p \geq D_1$ for all $p \in [1,\infty)$, it is enough to consider the case $p=1$. Let $x$, $y\in C$ and a representing coupling $\pi \in \Pc_{(x,y)}$. By definition of~\eqref{eqdef:lipSet} we have for any $a \in A_{d\text{-Lip}}(C)$ that
\begin{multline*}
\iint_{\Ext C \times \Ext C} d(z,z_*) \, \dd \pi(z,z_*) \geq \iint_{\Ext C \times \Ext C} |a(z)-a(z_*)| \, \dd \pi(z,z_*) \\
\geq \left| \iint_{\Ext C \times \Ext C} (a(z)-a(z_*)) \, \dd \pi(z,z_*)\right|.
\end{multline*}
Since $a$ is affine and continuous on $C$ and $\pi \in \Pc_{(x,y)}$, the right-most hand-side above reduces to
\[
\left| a(x) - a(y)\right|.
\]
Taking the infimum of the left-hand-side as $\pi$ runs through $\Pc_{(x,y)}$ yields
\begin{equation}
    \pRoof{1}(x,y) \geq |a(x)-a(y)|. 
\end{equation}
Since the above holds for all $(x,y) \in C\times C$, we can consider any chain, with $N\geq 1$, $ x^1 = x$, $x^n \in C$ and $x^{N+1} = y$, and obtain
\[
\sum_{n=1}^N \pRoof{1} (x^n,x^{n+1}) \geq \sum_{n=1}^N |a(x^n)-a(x^{n+1})| \geq \left| a(x) - a(y)\right|.
\]
Taking the infimum over all chains, we conclude that $D_1(x,y) \geq |a(x) - a(y)|$.

\medskip

\noindent Now to prove~\eqref{eq:lowerBdDp}, we simply notice that the hypothesis~\eqref{eq:lowerDbd} ensures that
\[
\left\{ f \in X^*,\; \; \|f\|_{X^*} = 1\right\} \subset A_{d\text{-Lip}}(C),
\]
and we conclude by using the duality formula of the norm
\[
\forall x,y \in C, \qquad \|x-y\| = \underset{\|f\|_{X^*}=1}{\sup_{f \in X^*}}f(x-y) = \underset{\|f\|_{X^*}=1}{\sup_{f \in X^*}}f(x) - f(y). \qquad \phantom{x}
\]
\end{proof}

\medskip

We conclude this section with the proof of Theorem~\ref{theo:metricDp}.

\begin{proof}[Proof of Theorem~\ref{theo:metricDp}]
Using Propositions~\ref{prop:propertiesDp} and~\ref{prop:sepLipschitz}, $D_p$ is a distance on $C$ under the assumption~\eqref{eq:sepHyp}. All conclusions of the theorem follow from Proposition~\ref{prop:propertiesDp} (basic properties), Proposition~\ref{prop:extensionproperty} (extension), Proposition~\ref{prop:monotonicity} (monotonicity in $p$), Proposition~\ref{prop:hatDpContinuousGeneral} (continuity with respect to the natural topology) and Proposition~\ref{prop:structureCDp} (metric structure). Finally~\eqref{eqTheo:ExtendsYes} comes by combining the sub-extension property of $D_p$ with the result~\eqref{eq:lowerBdDp} of Proposition~\ref{prop:sepLipschitz}.
\end{proof}

%% file: Section_applications.tex
\section{Applications to optimal transport on classical and quantum structures} \label{sec:applications}

\noindent This section is self-contained assuming the results of Section~\ref{sec:ChoquetWasserstein}, whose proofs are given in Section~\ref{sec:proofs}. We apply the construction to the classical state space of the simplex (Subsection~\ref{subsec:classical}), the semiclassical case of Golse--Paul, a product of quantum and classical state spaces (Subsection~\ref{subsec:GPsemiclassical}), and the entanglement-free quantum case (Subsection~\ref{sec:quantumChoquetWasserstein}), related to the construction of Beatty and Stilck Fran\c{c}a~\cite{beatty2025order} (see also T\'oth and Pitrik~\cite{toth2023quantum,toth2025quantum}). Throughout, for $x \in C$ whose resultant map is $r$, we denote by $\Pc_x := r^{-1}(\{x\})$ its set of representing measures (see Definition~\ref{def:representants}).

\subsection{Classical optimal transport: the case of the simplex} \label{subsec:classical}
Classical optimal transport focuses on the construction of costs, or distances, on spaces of probability measures. In the framework of Choquet theory~\cite{phelps2002lectures}, a compact metrizable convex set is a \emph{simplex} when its corresponding resultant map is \emph{injective}, meaning that every point in the set has a unique representation as a convex combination of extreme points. By Bauer~\cite{bauer1961silovscher}, any such simplex whose set of extreme points is closed is affinely homeomorphic to the set of Radon probability measures over its extreme points. In other words,
\[
C \text{ compact metrizable simplex, } {\rm Ext}\, C \text{ closed} \iff C = \Pc(E), \ E \text{ compact},
\]
in which case we have
 \[
 {\rm Ext} C =  \{\delta_x\}_{x \in E} \cong E.
 \]
Along with the surjectivity given by Choquet's theorem, $r$ is then a bijection between $\Pc({\rm Ext}\,C)$ and $C$, so that every point of $C$ has a unique representing measure, and therefore
\[
\pRoof{p} \; = \; D_p \; = \; W_p,
\]
the standard Wasserstein-$p$ distance on $\Pc({\rm Ext}\,C) \cong C$. Classical optimal transport thus amounts to performing convex extensions over a simplex, recovering the result of Savar\'e and Sodini~\cite{savare2022simple}, in which the standard Kantorovich problem is recast as extending a cost from the extreme boundary of a (Bauer) simplex to the whole simplex. 

\subsection{Semiclassical optimal transport \semiMetric of Golse and Paul} 
\label{subsec:GPsemiclassical} In the series of papers~\cite{golse2017schrodinger,golse2021semiclassical,golse2022quantitative}, Golse and Paul have developed and used an optimal transport-like cost between quantum density matrices and classical probability densities. This setting is not compact and we here only recast the Golse--Paul cost as a folded Kantorovich cost (Proposition~\ref{prop:GolsePaulRecast}). We recall the
\begin{definition}[Golse and Paul semiclassical \semiMetric~\cite{golse2017schrodinger}]
Let $\lambda, \hbar > 0$, define for all $(x,\xi) \in \R^d \times \R^d$ the cost operator $\hat{c}_{\lambda, \hbar}(x,\xi)$ over $L^2(\R^d)$ by, for any $\psi \in H^2(\R^d)$ and $y \in \R^d$,
\[
\hat{c}_{\lambda, \hbar}(x,\xi) \, \psi(y) := \lambda |x-y|^2 \psi(y) + \left[(\xi+i \hbar \nabla_y)^2 \psi \right](y).
\]
For a given trace one positive semidefinite self-adjoint operator $R$ over $L^2(\R^d)$ and a probability density (w.r.t. the Lebesgue measure) $f$ over $\R^d \times \R^d$, we say that a map $(x,\xi) \in \R^d \times \R^d \mapsto Q(x,\xi)$ is a coupling between $R$ and $f$ if $Q(x,\xi)$ is a trace-class positive semidefinite self-adjoint operator over $L^2(\R^d)$ for all $(x,\xi)$, and satisfies the marginal constraints
\[
\mathrm{Tr}(Q(x,\xi)) = f(x,\xi), \qquad \iint_{\R^d \times \R^d} Q(x,\xi) \, \dd x \, \dd \xi = R. 
\]
The Golse and Paul semiclassical \semiMetric between $R$ and $f$ is then defined by
\begin{equation}
    \mathrm{GP}_{\lambda, \hbar}(R,f) := \inf_{Q \emph{ coupling between } R \emph{ and } f} \; \; \left(\iint_{\R^d \times \R^d} \mathrm{Tr} \left(\hat{c}_{\lambda,\hbar}(x,\xi) \, Q(x,\xi) \right) \, \dd x \, \dd \xi \right)^\frac12.
\end{equation}
\end{definition}
\noindent We can recast the Golse and Paul \semiMetric as a maximal convex extension (folded Kantorovich cost).
\begin{proposition}\label{prop:GolsePaulRecast} Consider the convex sets $C_1 = \mathfrak{S}_{L^2(\R^d)}$ of trace one positive semidefinite self-adjoint operators over $L^2(\R^d)$, and $C_2 = \Pc(\R^d \times \R^d)$, whose extreme boundaries are respectively $E_1 = \mathbf{P}_{L^2(\R^d)}$ (the rank-one orthogonal projectors over $L^2(\R^d)$), and $E_2 \cong \R^d \times \R^d$. For $\lambda, \hbar > 0$, define the cost function
\begin{equation}
c_{\lambda, \hbar} : \left( P_\psi, \, (x,\xi) \right) \in E_1 \times E_2 \; \mapsto \; \langle \psi \; | \; \hat{c}_{\lambda, \hbar}(x,\xi) \; | \; \psi \rangle_{L^2(\R^d)} \in \R_+ \cup \{+\infty\}.
\end{equation}
Then for any $R \in \mathfrak{S}_{L^2(\R^d)}$ and probability density $f$ on $\R^d \times \R^d$, we have, denoting by $\nu_f$ the measure $f \, \dd x \, \dd \xi$, 
    \begin{equation} \label{eq:GPrecast}
        \mathrm{GP}_{\lambda, \hbar}(R,f)^2 = (c_{\lambda, \hbar})_*(R,\nu_f),
    \end{equation}
where $(c_{\lambda,\hbar})_*$ is the maximal convex extension of $c_{\lambda,\hbar}$ from $E_1 \times E_2$ to $C_1 \times C_2$, see~\eqref{eq:foldedKantoRelation}.
\end{proposition}
As a consequence of~\eqref{eq:GPrecast}, $\mathrm{GP}_{\lambda,\hbar}$ naturally extends to any Radon probability measures $\nu \in \Pc(\R^d \times \R^d)$, instead of sole densities $f$.

\begin{proof}
    Let $\mu \in \Pc(\mathbf{P}_{L^2(\R^d)})$ represent $R$ and $\pi$ a coupling between $\mu$ and $\nu_f$. We have
    \begin{align*}
      \iint_{E_1 \times E_2} c_{\lambda, \hbar}(e_1,e_2) \, \dd \pi(e_1,e_2) &= \iint_{\mathbf{P}_{L^2(\R^d)} \times \R^{2d}} \langle \psi \; | \; \hat{c}_{\lambda, \hbar}(x,\xi) \; | \; \psi \rangle_{L^2(\R^d)} \, \dd \pi(P_{\psi},(x,\xi)) \\
      &= \iint_{\mathbf{P}_{L^2(\R^d)} \times \R^{2d}} \mathrm{Tr} \left(\hat{c}_{\lambda, \hbar}(x,\xi) \; P_\psi \right) \, \dd \pi(P_{\psi},(x,\xi)).
    \end{align*}
By the disintegration theorem, since $\nu_f$ is the second marginal of $\pi$, there exists a family of probability measures over $\mathbf{P}_{L^2(\R^d)}$ indexed by $(x,\xi)$ and denoted $\pi(\cdot |x,\xi)$, such that
\[
\dd \pi(P_{\psi},(x,\xi)) = \dd \pi(P_\psi |x,\xi) \; f(x,\xi) \, \dd x \, \dd \xi.
\]
Therefore the last integral becomes
\begin{equation*}
\iint_{\mathbf{P}_{L^2(\R^d)} \times (\R^d \times \R^d)} \mathrm{Tr} \left(\hat{c}_{\lambda, \hbar}(x,\xi) \; P_\psi \right) \, \dd \pi(P_\psi |x,\xi) \, f(x,\xi) \, \dd x \, \dd \xi,
\end{equation*}
that is
\begin{equation*}
\iint_{\R^d \times \R^d} \mathrm{Tr} \left(\hat{c}_{\lambda, \hbar}(x,\xi) \; Q(x,\xi) \right) \,   \dd x \, \dd \xi,
\end{equation*}
with
\[
Q(x,\xi) :=  f(x,\xi)\int_{\mathbf{P}_{L^2(\R^d)}} P_\psi \; \dd \pi(P_\psi |x,\xi).
\]
As $\pi(\cdot |x,\xi)$ is a probability measure over $\mathbf{P}_{L^2(\R^d)}$, and since the first marginal of $\pi$ is $\mu$ representing~$R$, $Q$ is a valid coupling between $R$ and $f$. In turn, we have
\[
\iint_{E_1 \times E_2} c_{\lambda, \hbar}(e_1,e_2) \, \dd \pi(e_1,e_2) \geq \mathrm{GP}_{\lambda, \hbar}(R,f)^2,
\]
and taking the infimum over all coupling measures $\pi$ and all representing probabilities $\mu$ yields
\[
\inf_{\mu \in \Pc_R} \; \mathcal{K}_{c_{\lambda, \hbar}}(\mu, \nu_f) \geq \mathrm{GP}_{\lambda, \hbar}(R,f)^2.
\]
Conversely, consider some coupling $Q$ between $R$ and $f$. For any $(x,\xi) \in \R^d \times \R^d$, since $Q(x,\xi)$ is a positive semidefinite self-adjoint operator over $L^2(\R^d)$ with finite trace $f(x,\xi)$, there exists (by the spectral theorem)%\footnote{The decomposition may be chosen measurably in $(x,\xi)$, by standard measurable-selection theorems, see e.g.~\cite{castaing1977convex}}
a Radon probability measure $\dd \pi(\cdot |x,\xi)$ over $\mathbf{P}_{L^2(\R^d)}$ such that
\[
Q(x,\xi) =  f(x,\xi)\int_{\mathbf{P}_{L^2(\R^d)}} P_\psi \; \dd \pi(P_\psi |x,\xi).
\]
We let
\[
\dd \pi (P_\psi, (x,\xi)) := \dd \pi(P_\psi |x,\xi) \, f(x,\xi) \, \dd x \, \dd \xi, \quad \dd \mu(P_\psi) := \iint_{(x,\xi) \in \R^d \times \R^d}\dd \pi(P_\psi,(x,\xi)).
\]
Then $\mu$ is a Radon probability measure over $\mathbf{P}_{L^2(\R^d)}$, $\pi$ is a coupling between $\mu$ and $\nu_f$, and it is easy to check with a Fubini argument that $\mu$ represents $R$. Then, as again
\[
\iint_{E_1 \times E_2} c_{\lambda, \hbar}(e_1,e_2) \, \dd \pi(e_1,e_2) = \iint_{\R^d \times \R^d} \mathrm{Tr} \left(\hat{c}_{\lambda, \hbar}(x,\xi) \; Q(x,\xi) \right) \,   \dd x \, \dd \xi,
\]
we obtain
\[
\iint_{\R^d \times \R^d} \mathrm{Tr} \left(\hat{c}_{\lambda, \hbar}(x,\xi) \; Q(x,\xi) \right) \,   \dd x \, \dd \xi \geq  \inf_{\mu \in \Pc_R} \; \mathcal{K}_{c_{\lambda, \hbar}}(\mu, \nu_f).
\]
As this holds for any coupling $Q$, taking the infimum over $Q$ on the left-hand-side yields
\[
\mathrm{GP}_{\lambda, \hbar}(R,f)^2 \geq \inf_{\mu \in \Pc_R} \; \mathcal{K}_{c_{\lambda, \hbar}}(\mu, \nu_f),
\]
which concludes the proof.
\end{proof}

\subsection{Entanglement-free quantum optimal transport of Beatty and Stilck Fran\c{c}a} \label{sec:quantumChoquetWasserstein}
\noindent Let us now specialize our construction to the convex set of density matrices, the quantum state space.
\subsubsection{Preliminary definitions} \label{subsec:densityMatricesIntro} Let $\hilbert$ be a separable complex Hilbert space endowed with a scalar product $\braket{\cdot | \cdot}$, and $\B(\hilbert)$ the set of bounded operators on $\hilbert$.

\begin{definition}[Density matrices and pure states] \label{def:densityMatricesandPureStates}
We consider the set of \emph{density matrices} $\dens$ to be the set of trace $1$, positive semidefinite self-adjoint operators over $\hilbert$, i.e.
\begin{equation} \label{eqdef:densitymatrices}
\dens := \left\{ \rho \in \B (\hilbert) \; \; | \; \; \rho  = \rho^\dagger, \; \; \forall \psi \in \hilbert, \, \braket{\psi \, | \, \rho \, | \,\psi}  \geq 0 \; \text{ and } \; \mathrm{Tr}(\rho) = 1\right\}.
\end{equation}
The set of \emph{pure states} $\purest$ (identified with the complex projective space of $\hilbert$) is the set of rank-one orthogonal projectors over $\hilbert$, i.e.
\begin{equation} \label{eqdef:purestates}
\purest := \left\{ P_\psi : \phi \in \hilbert \mapsto \braket{\psi | \phi} \; \psi \; \; | \; \; \psi \in \hilbert \text{ and } \langle\psi | \psi \rangle = 1 \right\}.
\end{equation}
\end{definition}

\noindent The set $\dens$ is a convex subset of $\B (\hilbert)$, whose extreme boundary is~$\purest$. For any $\rho \in \dens$, the spectral theorem ensures the existence of a family $\{\lambda_n\}_n \subset \R_+$ and an orthonormal basis $\{\psi_n\}_n $ of $\hilbert$ such that $\sum_n \lambda_n = 1$ and
\[
\rho = \sum_{n} \lambda_n \, P_{\psi_n}.
\]
Then the measure $\displaystyle \sum_n \lambda_{n} \, \bm{\delta}_{P_{\psi_n}}$, with $\bm{\delta}_P$ denoting the Dirac mass at $P \in \purest$, is a probability measure over $\purest$ which represents $\rho$. In turn, the spectral theorem may be seen as a stronger version of Choquet's Theorem in the spectral setting. Denoting by $r : \Pc(\purest) \to \dens$ the resultant map, which sends a representing measure to its barycenter, we recover the Choquet identification of Section~\ref{sec:ChoquetWasserstein},
\begin{equation} \label{eq:ChoquetlRepresentationQuantum}
    \dens \; \; \cong \; \;\faktor{\Pc(\purest)}{r}.
\end{equation}

\subsubsection{Quantum folded Wasserstein} \label{subsec:quantumChoWass}
Throughout this subsection, $\HH$ denotes a complex Hilbert space of finite dimension $\cong \mathbb{C}^n$, and the space $\B(\HH)$ of bounded operators is endowed with some norm $\|\cdot\|$ whose topology we call the natural topology (all norms are equivalent in finite dimension). We consider a distance $d$ on the set of pure states $\purest$ continuous with respect to the natural topology (hence inducing it, $\purest$ being compact, see Section~\ref{sec:proofs}). The set $\dens$ is then a compact (since $\dim\HH<\infty$) and stable (see Section~\ref{sec:ChoquetWasserstein}) convex subset of $\B(\HH)$, whose extreme boundary is $\purest$. Our construction of Section~\ref{sec:ChoquetWasserstein} therefore applies, with $X=\B(\HH)$, $C=\dens$ and $\Ext C=\purest$.

\medskip

\noindent For $p \in [1,\infty)$, we denote by $W_p$ the standard Wasserstein-$p$ distance on $\Pc(\purest)$ built from $(\purest,d)$. 
\begin{definition}[Quantum folded Wasserstein] \label{def:quantumChoquetWass} The quantum folded Wasserstein\nobreakdash-$p$~\psDistance on $\dens$ associated to $d$ (see Section~\ref{sec:ChoquetWasserstein}, Equations~\eqref{eqdef:PhiMaxConvExt}--\eqref{eqdef:PhimaxMetricconvexExtension}), writes in this setting, for any $\rho,\sigma \in \dens$,
\begin{equation} \label{eqdef:Dpquantum}
D_p (\rho, \sigma) = \inf \left\{ \sum_{n=1}^N\pRoof{p} (\gamma_n,\gamma_{n+1}), \quad N\geq 1, \quad  \gamma_1 = \rho, \quad \gamma_n \in \dens, \quad \gamma_{N+1} = \sigma \right\},
\end{equation}
where $\pRoof{p}$ is the quantum folded Kantorovich-$p$ \semiDistance
\begin{equation} \label{eqdef:Dphatquantum}
\pRoof{p}(\rho, \sigma) = \inf_{\mu \in \Pc_{\rho}} \; \; \; \inf_{\nu \in \Pc_{\sigma}} W_p(\mu,\nu),
\end{equation}
with $\Pc_{\rho} \subset \Pc(\purest)$ the set of representing probabilities of $\rho$ (see Definition~\ref{def:representants}) defined by
\[
\Pc_\rho := \left\{\mu \in \Pc(\purest) \;  | \; \; \int_{\purest} P \, \dd \mu(P) = \rho \right\}.
\]
We also recall the definition of the set of representing couplings of $(\rho,\sigma) \in \dens \times \dens$ (introduced at the beginning of Section~\ref{sec:proofs}),
\[
\Pc_{(\rho,\sigma)} := \left\{\pi \in \Pc(\purest \times \purest) : \int_{\purest \times \purest} P \, \dd \pi(P,Q) = \rho, \ \int_{\purest \times \purest} Q \, \dd \pi(P,Q) = \sigma \right\}.
\]
\end{definition}

\noindent We recall Proposition~\ref{prop:propertiesDp} in Section~\ref{sec:ChoquetWasserstein}, ensuring in particular that $D_p$ is indeed a pseudo-distance, $\pRoof{p}$ is indeed a semi-distance\footnote{We recall that a pseudo-distance may fail to separate points, while a semi-distance may fail to satisfy the triangle inequality.}, $D_p \leq \pRoof{p}$ and that $\pRoof{p}$ is a distance if and only if $\pRoof{p} = D_p$. We highlight that $\pRoof{p}$ is exactly the Beatty and Stilck Fran\c{c}a \semiDistance~\cite{beatty2025order}, as discussed later in Subsection~\ref{subsec:Beatty-Franca}. It is the $p$-th root of the convex extension of $d^p$ on $C \times C$, i.e. $\pRoof{p}^p=(d^p)_*$. This is the same mechanism as for entanglement measures, which are convex roofs of functions defined on pure states~\cite{uhlmann2010roofs,vollbrecht2001entanglement,chitambar2019quantum}, here applied to a distance.

\medskip

We now state our main Theorem~\ref{theo:mainquantum} of this section, which follows from Theorem~\ref{theo:metricDp} and the fact that $\dens$ is a compact and stable convex set. Since the norm $\|\cdot\|$ can be chosen arbitrarily, the required assumption~\eqref{eq:lowerDbdquantum} is satisfied by any norm on $\B(\hilbert)$ as well as by the Fubini-Study distance (take $\|\cdot\|$ to be $2^{-\frac12}$ times the Frobenius norm).

\begin{theorem} \label{theo:mainquantum}
Let $p \in [1,\infty)$ and further assume $d$ to dominate the norm,
\begin{equation} \label{eq:lowerDbdquantum}
\|P-Q\| \leq d(P,Q) \qquad \text{for all } P,Q \in \purest.
\end{equation}
Then $D_p$ is a distance on $\dens$ whose induced topology is equivalent to the natural topology; in particular $(\dens, D_p)$ is a compact Polish space.

\smallskip

\noindent Moreover, $(\dens,D_1)$ is always geodesic, and if $(\purest, d)$ is geodesic, then $(\dens, D_p)$ is also geodesic for all $p>1$.

\smallskip

\noindent Finally, $D_p$ dominates the norm, $D_p(\rho,\sigma) \geq \|\rho-\sigma\|$ for all $\rho,\sigma \in \dens$, and sub-extends $d$, ${D_p}_{|\purest \times \purest} \leq d$, with equality when $d(P,Q)=\|P-Q\|$ on $\purest$.
\end{theorem}

As examples, we provide in the next corollary the properties of the space $(\dens,D_p)$ in the case of two typical distances over $\purest$.

\begin{corollary}\label{cor:quantumExamples}
When the distance over $\purest$ is the \emph{Frobenius} distance
\begin{equation} \label{eq:dFr}
d^{\mathrm{Fr}}(P_{\psi}, P_{\phi}) = \sqrt{2}\,\sqrt{1 -|\braket{\psi|\phi}|^2} = \|P_\psi - P_\phi\|_{\mathrm{Fr}},
\end{equation}
$D^{\mathrm{Fr}}_p$ is a distance on $\dens$ whose induced topology is equivalent to the natural topology, \emph{extends} $d^{\mathrm{Fr}}$, that is, ${D^{\mathrm{Fr}}_p}_{|\purest \times \purest} = d^{\mathrm{Fr}}$, and dominates the Frobenius norm. Moreover $(\dens, D^{\mathrm{Fr}}_1)$ is geodesic.

\medskip

\noindent When it is the \emph{Fubini--Study} distance
\begin{equation} \label{eq:dFS}
d^{\mathrm{FS}}(P_{\psi}, P_{\phi}) = \arccos |\braket{\psi|\phi}|,
\end{equation}
$D^{\mathrm{FS}}_p$ is a distance on $\dens$ whose induced topology is equivalent to the natural topology, \emph{sub-extends} $d^{\mathrm{FS}}$, that is, ${D^{\mathrm{FS}}_p}_{|\purest \times \purest} \leq d^{\mathrm{FS}}$, and dominates $2^{-1/2}\|\cdot\|_{\mathrm{Fr}}$. Moreover $(\dens, D^{\mathrm{FS}}_p)$ is geodesic for every $p\geq1$.
\end{corollary}

\begin{proof}
By definition of the Frobenius norm $\|\cdot\|_{\text{Fr}}$ on $\B(\hilbert)$, we have for any $P_\psi \in \purest$ that
\[
\|P_\psi-P_\phi\|_{\text{Fr}}^2 = \mathrm{Tr} \left((P_\psi-P_\phi) (P_\psi-P_\phi)^\dagger \right)  = 2 - 2 \mathrm{Tr}(P_\psi P_\phi) = 2 (1 -  |\braket{\psi|\phi}|^2),
\]
and as such,
\[
d^{\text{Fr}}(P_\psi, \, P_\phi) = \|P_\psi-P_\phi\|_{\text{Fr}}.
\]
Therefore $d^{\text{Fr}}$ obviously induces the natural topology and upper-bounds $\|\cdot\|_{\text{Fr}}$. Applying Theorem~\ref{theo:mainquantum} with the choice $\|\cdot\|=\|\cdot\|_{\mathrm{Fr}}$ in~\eqref{eq:lowerDbdquantum} then gives that $D^{\mathrm{Fr}}_p$ is a distance whose induced topology is equivalent to the natural topology, dominates the Frobenius norm, and sub-extends $d^{\mathrm{Fr}}$ with equality on $\purest$, since $d^{\mathrm{Fr}}$ is the chosen norm distance itself. The same theorem yields that $(\dens,D^{\mathrm{Fr}}_1)$ is geodesic.

\medskip

\noindent Turning to the Fubini--Study distance, note that
\[
\sin \left(d^{\text{FS}}(P_\psi,P_\phi) \right) = \sqrt{1 - |\braket{\psi|\phi}|^2} = \frac{1}{\sqrt{2}} \, \|P_\psi-P_\phi\|_{\text{Fr}},
\]
implying
\[
 \frac{1}{\sqrt{2}} \, \|P_\psi-P_\phi\|_{\text{Fr}}\leq d^{\text{FS}}(P_\psi,P_\phi) \leq  \frac{\pi}{2\sqrt{2}} \, \|P_\psi-P_\phi\|_{\text{Fr}}.
\]
Therefore $d^{\text{FS}}$ induces the natural topology and upper-bounds $2^{-\frac12}\|\cdot\|_{\text{Fr}}$. Applying Theorem~\ref{theo:mainquantum} with the choice $\|\cdot\|=2^{-\frac12}\|\cdot\|_{\mathrm{Fr}}$ in~\eqref{eq:lowerDbdquantum} gives that $D^{\mathrm{FS}}_p$ is a distance whose induced topology is equivalent to the natural topology, dominates $2^{-\frac12}\|\cdot\|_{\mathrm{Fr}}$, and sub-extends $d^{\mathrm{FS}}$. As $(\purest,d^{\text{FS}})$ is moreover geodesic, the geodesy statement of the theorem yields that $(\dens, D^{\mathrm{FS}}_p)$ is geodesic for $p>1$, while $(\dens, D^{\mathrm{FS}}_1)$ is geodesic in any case.
\end{proof}

\noindent For $p=2$, the Frobenius case can be sharpened: $D_2^{\mathrm{Fr}}$ actually coincides with $\pRoof{2}^{\mathrm{Fr}}$~\cite{miller2026metricity}, see Remark~\ref{rem:subadditivityBSF}.

\subsubsection{Beatty and Stilck Fran\c{c}a\texorpdfstring{~\cite{beatty2025order}}{} and quantum folded Wasserstein}
\label{subsec:Beatty-Franca} 
We recall that we have assumed here $\HH$ finite-dimensional and $d$ is continuous on $\purest\times\purest$ with respect to the natural topology. Given $\rho, \sigma \in \dens$, the set of Beatty and Stilck Fran\c{c}a transport plans between $\rho$ and $\sigma$, here denoted by $\mathcal{C}^{\mathrm{BSF}}(\rho, \sigma)$, is defined\footnote{In~\cite{beatty2025order}, unit vectors $\psi_j$ are used instead of pure states $P_{\psi_j}$, but it is equivalent here because $d$ is defined on pure states.} as
\begin{multline}\label{eqdef:BFTransportPlans}
\mathcal{C}^{\mathrm{BSF}}(\rho, \sigma) := \Big\{ (q_j, P_{\psi_j}, P_{\psi^*_j})_{1 \leq j \leq J} \in (\R_+^* \times \purest \times \purest)^J \; : \; J \in \N^*, \\
\sum_{j=1}^{J} q_j \, P_{\psi_j} = \rho, \ \sum_{j=1}^{J} q_j \, P_{\psi^*_j} = \sigma \Big\}.
\end{multline}
Then the $p$-Beatty and Stilck Fran\c{c}a \semiDistance is defined in~\cite{beatty2025order} by
\begin{equation} \label{eqdef:BFdist}
D^{\mathrm{BSF}}_p(\rho, \sigma) := \left( \inf_{(q_j, P_{\psi_j}, P_{\psi^*_j})_{1 \leq j \leq J} \in \mathcal{C}^{\mathrm{BSF}}(\rho, \sigma)} \, \sum_{j=1}^{J} q_j \, d(P_{\psi_j}, P_{\psi^*_j})^p \right)^{\frac{1}{p}}.
\end{equation}
The link between our construction and the one of Beatty and Stilck Fran\c{c}a lies in the following proposition.
\begin{proposition} \label{prop:BFrecast}
For every $p \in [1,\infty)$, the order-$p$ Beatty--Stilck Fran\c{c}a semi-distance~\eqref{eqdef:BFdist} coincides with the folded Kantorovich-$p$ semi-distance:
\begin{equation} \label{eq:Beatty}
D^{\mathrm{BSF}}_p = \pRoof{p}.
\end{equation}
\end{proposition}
\begin{proof}
Recall that Proposition~\ref{prop:linearProgramhatDp} ensures that $\pRoof{p}^p$ is the convex extension on $\dens \times \dens$ of $d^p$ in the sense of 
\begin{equation*}
\pRoof{p}(\rho, \sigma)^p = \inf_{\pi \in \Pc_{(\rho,\sigma)}} \iint_{\purest \times \purest} d(P,Q)^p \, \dd \pi(P,Q).
\end{equation*}
We can identify
\[
(q_j, P_{\psi_j}, P_{\psi^*_j})_{1 \leq j \leq J} \in (\R_+^* \times \purest \times \purest)^J   \qquad \text{with} \qquad \sum_{j=1}^{J} q_j \, \bm{\delta}_{P_{\psi_j} \otimes P_{\psi_j^*}} \in \Pc(\purest \times \purest),
\]
so that
\[
\mathcal{C}^{\mathrm{BSF}}(\rho, \sigma) \cong \left\{\pi \in \Pc_{(\rho,\sigma)}  \; \; |  \; \; \pi\text{ is finitely supported} \right\}
\]
and~\eqref{eqdef:BFdist} can be recast as
\[
D^{\mathrm{BSF}}_p(\rho, \sigma)^p = \underset{\pi\text{ is finitely supported}}{\inf_{\pi \in \Pc_{(\rho,\sigma)}}} \, \iint_{\purest \times \purest} d(P, Q)^p \, \dd \pi(P,Q).
\]
The equality of the two quantities follows from the existence of a finitely supported optimal representing coupling, Proposition~\ref{prop:quantumoptimalplan} below (a consequence of Proposition~\ref{prop:linearProgramhatDp}).
\end{proof}

\begin{remark}
A representing coupling $\pi \in \Pc_{(\rho,\sigma)}$ represents a \emph{separable bipartite} coupling between $\rho$ and $\sigma$
\[
\iint_{\purest \times \purest} P \otimes Q \, \dd \pi(P,Q).
\]
Entangled bipartite states are therefore not considered in our construction.
\end{remark}

\begin{remark}[On the subadditivity of $D^{\mathrm{BSF}}_p$]\label{rem:subadditivityBSF}
Identifying $D^{\mathrm{BSF}}_p$ with $\pRoof{p}$ provides information on the triangle inequality. The double inf formulation~\eqref{eqdef:Dphatquantum} gives no reason for $\pRoof{p}$ to be subadditive in general, and by Proposition~\ref{prop:propertiesDp},
\[
\pRoof{p} \equiv D^{\mathrm{BSF}}_p \text{ is a distance} \quad \iff \quad \pRoof{p} = D_p .
\]
Both objects have their advantages and drawbacks. On the one hand, $\pRoof{p} \equiv D^{\mathrm{BSF}}_p$ is computable (Proposition~\ref{prop:quantumoptimalplan}) and always extends $d$, but may fail to be subadditive. On the other hand, $D_p$ is, in the cases of Theorem~\ref{theo:mainquantum} and Corollary~\ref{cor:quantumExamples}, a genuine distance with good metric properties, but is harder to compute and, outside the norm case, a priori only sub-extends $d$. The optimal situation is $D_p = \pRoof{p} \equiv D^{\mathrm{BSF}}_p$, where all desirable properties  hold at once: a computable distance that extends $d$, with no need for  the lower bound~\eqref{eq:lowerDbdquantum}. This holds in particular under the criterion of Proposition~\ref{prop:subadditivitygeneral}.
This optimal situation is known to hold for $p=2$ and the Frobenius distance $d^{\mathrm{Fr}}$ of Corollary~\ref{cor:quantumExamples}: Miller~\cite{miller2026metricity} recently proved, in every finite dimension, that $\pRoof{2}^{\mathrm{Fr}} \equiv D^{\mathrm{BSF}}_2$ satisfies the triangle inequality, hence $\pRoof{2}^{\mathrm{Fr}} = D_2^{\mathrm{Fr}}$. The proof relies on a dual formulation of the convex roof, and the cases $p \neq 2$ or other ground distances remain open.
\end{remark}

\paragraph{Related results on $D^{\mathrm{BSF}}_p \equiv \pRoof{p}$}
Our setting allows one to recover a few results of~\cite{beatty2025order}: in the following proposition, we almost recover the result of~\cite[Proposition~16]{beatty2025order}, passing from an upper-bound of $2 (\dim \hilbert)^2$ to $2(\dim \hilbert)^2 - 1$ atoms for the optimal transport plan.
\begin{proposition} \label{prop:quantumoptimalplan}
Let $p \in [1,\infty)$. For all $\rho, \sigma \in \dens$, $\pRoof{p}(\rho, \sigma)$ is realized by a finitely supported optimal representing coupling $\pi^* \in \Pc_{(\rho,\sigma)}$ with at most $2(\dim \hilbert)^2 - 1$ atoms.
\end{proposition}

\begin{proof}
Take $\rho \in \dens$ and denote by $(\lambda_n)_{1 \leq n \leq \dim \hilbert} \in (\R_+)^{\dim \hilbert}$ its eigenvalues, along with an associated orthonormal eigenbasis $(\psi_n)_{1 \leq n \leq \dim \hilbert}$. Let us set, for all $1 \leq n,m \leq \dim \hilbert$,
\[
g_n : P_\phi \in \purest \quad \mapsto \quad \langle \psi_n | P_{\phi} | \psi_n \rangle  \in \R, \qquad f_{n,m} : P_\phi \in \purest \quad \mapsto \quad \langle \psi_n | P_{\phi} | \psi_m \rangle  \in \mathbb{C},
\]
which are continuous and bounded respectively real-and complex-valued functions over $\purest$. Note that the set $\Pc_{\rho}$ of representing measures of $\rho$ (see Definition~\ref{def:representants}) can be recast as
\begin{multline}
\Pc_{\rho} = \Big\{ \mu \in \Pc(\purest) \; : \; \int_{\purest} g_{n} \, \dd \mu = \lambda_n \text{ for } 1 \leq n \leq \dim \hilbert - 1, \\
\int_{\purest} f_{n,m} \, \dd \mu =0 \text{ for } 1 \leq n < m\leq \dim \hilbert \Big\}.
\end{multline}
Since $\mu$ sums to one, it is enough to consider only $\dim \HH - 1$ elements of the family $(g_n)_{1 \leq n \leq \dim \hilbert}$. The number of \emph{real} linear constraints is therefore $\dim \HH - 1 + 2\frac{\dim \HH(\dim \HH -1)}{2} = (\dim \HH)^2 - 1$. Then the proposition comes as a corollary of Proposition~\ref{prop:linearProgramhatDp}.
\end{proof}
Second, we answer by the negative a question raised by the authors in the proof of~\cite[Corollary~13]{beatty2025order}, asking whether H\"older continuity of the $2$-norm with respect to $d$ is necessary for $D^{\mathrm{BSF}}_p = \pRoof{p}$ to separate points, with the following proposition (following from Proposition~\ref{prop:propertiesDp} in the general setting).
\begin{proposition} Let $p \in [1,\infty)$. Then $\pRoof{p} \equiv D^{\mathrm{BSF}}_p$ separates points on $\dens$.
\end{proposition}
Lastly, we also recover the uniform continuity\footnote{Since $\dens$ is compact in the natural topology, continuity implies uniform continuity on $\dens$.} of $\pRoof{p}$ shown in~\cite[Proposition~17]{beatty2025order}, but with different tools which rely on the sole fact that the convex set $\dens$ is stable, see Proposition~\ref{prop:hatDpContinuousGeneral}.
\begin{proposition} Let $p \in [1,\infty)$. Then $\pRoof{p} \equiv D^{\mathrm{BSF}}_p$ is continuous on $\dens$ for the natural topology.  
\end{proposition}

%% file: biblio.bib
@article{papadopoulou1977geometry,
  title={On the geometry of stable compact convex sets},
  author={Papadopoulou, S.},
  journal={Math. Ann.},
  volume={229},
  number={3},
  pages={193--200},
  year={1977}
}

@article{corson1970metrizability,
  title={Metrizability of compact convex sets},
  author={Corson, H. H.},
  journal={Trans. Am. Math. Soc.},
  volume={151},
  number={2},
  pages={589--596},
  year={1970}
}

@article{volle2015duality,
  title={Duality for closed convex functions and evenly convex functions},
  author={Volle, M. and Martinez-Legaz, J. E. and Vicente-P{\'e}rez, J.},
  journal={J. Optim. Theory Appl.},
  volume={167},
  number={3},
  pages={985--997},
  year={2015}
}

@article{o1976openness,
  title={On the openness of the barycentre map},
  author={O'Brien, R. C.},
  journal={Math. Ann.},
  volume={223},
  number={3},
  pages={207--212},
  year={1976}
}

@article{shirokov2007strong,
  title={On the strong {CE}-property of convex sets},
  author={Shirokov, M. E.},
  journal={Math. Notes},
  volume={82},
  number={3},
  pages={395--409},
  year={2007}
}

@article{shirokov2012stability,
  title={Stability of convex sets and applications},
  author={Shirokov, M. E.},
  journal={Izv. Math.},
  volume={76},
  number={4},
  pages={840--856},
  year={2012}
}

@book{alfsen2012compact,
  title={Compact convex sets and boundary integrals},
  author={Alfsen, E. M.},
  volume={57},
  year={2012},
  publisher={Springer}
}

@article{winkler1988extreme,
  title={Extreme points of moment sets},
  author={Winkler, G.},
  journal={Math. Oper. Res.},
  volume={13},
  number={4},
  pages={581--587},
  year={1988}
}

@book{burago2001course,
  title={A course in metric geometry},
  author={Burago, D. and Burago, Y. and Ivanov, S.},
  volume={33},
  year={2001},
  publisher={American Mathematical Society Providence}
}

@article{golse2017schrodinger,
  title={The {S}chr{\"o}dinger equation in the mean-field and semiclassical regime},
  author={Golse, F. and Paul, T.},
  journal={Arch. Ration. Mech. Anal.},
  volume={223},
  pages={57--94},
  year={2017}
}

@article{golse2021semiclassical,
  title={Semiclassical evolution with low regularity},
  author={Golse, F. and Paul, T.},
  journal={J. Math. Pures Appl.},
  volume={151},
  pages={257--311},
  year={2021}
}

@article{golse2022quantitative,
  title={Quantitative observability for the {S}chr{\"o}dinger and {H}eisenberg equations: An optimal transport approach},
  author={Golse, F. and Paul, T.},
  journal={Math. Models Methods Appl. Sci.},
  volume={32},
  number={05},
  pages={941--963},
  year={2022}
}

@book{phelps2002lectures,
  title={Lectures on {C}hoquet's theorem},
  author={Phelps, R. R.},
  year={2002},
  publisher={Springer}
}

@article{bauer1961silovscher,
  title={{\v{S}}ilovscher {R}and und {D}irichletsches {P}roblem},
  author={Bauer, H.},
  journal={Ann. Inst. Fourier},
  volume={11},
  pages={89--136},
  year={1961}
}

@misc{peifer2020choquet,
  title={{C}hoquet theory},
  author={Peifer, C. J.},
  year={2020},
  howpublished={Bachelor's thesis, University of Copenhagen}
}

@article{trevisan2025quantum,
  title={Quantum optimal transport: an invitation},
  author={Trevisan, D.},
  journal={Boll. Unione Mat. Ital.},
  volume={18},
  number={1},
  pages={347--360},
  year={2025}
}

@article{carlen2014analog,
  title={An analog of the 2-{W}asserstein metric in non-commutative probability under which the fermionic {F}okker--{P}lanck equation is gradient flow for the entropy},
  author={Carlen, E. A. and Maas, J.},
  journal={Commun. Math. Phys.},
  volume={331},
  number={3},
  pages={887--926},
  year={2014}
}

@article{de2021quantum,
  title={Quantum optimal transport with quantum channels},
  author={De Palma, G. and Trevisan, D.},
  journal={Ann. Henri Poincar{\'e}},
  volume={22},
  pages={3199--3234},
  year={2021}
}

@article{golse2016mean,
  title={On the mean field and classical limits of quantum mechanics},
  author={Golse, F. and Mouhot, C. and Paul, T.},
  journal={Commun. Math. Phys.},
  volume={343},
  number={1},
  pages={165--205},
  year={2016}
}

@article{toth2023quantum,
  title={Quantum {W}asserstein distance based on an optimization over separable states},
  author={T{\'o}th, G. and Pitrik, J.},
  journal={Quantum},
  volume={7},
  pages={1143},
  year={2023}
}

@article{yan2012extension,
  title={Extension of convex function},
  author={Yan, M.},
  journal={J. Convex Anal.},
  volume={21},
  number={4},
  pages={965--987},
  year={2014}
}

@article{savare2022simple,
  title={A simple relaxation approach to duality for optimal transport problems in completely regular spaces},
  author={Savar{\'e}, G. and Sodini, G. E.},
  journal={J. Convex Anal.},
  volume={29},
  number={1},
  pages={1--12},
  year={2022}
}

@article{delon2020wasserstein,
  title={A {W}asserstein-type distance in the space of {G}aussian mixture models},
  author={Delon, J. and Desolneux, A.},
  journal={SIAM J. Imaging Sci.},
  volume={13},
  number={2},
  pages={936--970},
  year={2020}
}

@article{chitambar2019quantum,
  title={Quantum resource theories},
  author={Chitambar, E. and Gour, G.},
  journal={Rev. Mod. Phys.},
  volume={91},
  number={2},
  pages={025001},
  year={2019},
}

@article{uhlmann2010roofs,
  title={Roofs and convexity},
  author={Uhlmann, A.},
  journal={Entropy},
  volume={12},
  number={7},
  pages={1799--1832},
  year={2010}
}

@article{bucicovschi2010continuity,
  title={On the continuity and regularity of convex extensions},
  author={Bucicovschi, O. and Lebl, J.},
  journal={J. Convex Anal.},
  volume={20},
  number={4},
  pages={1113--1126},
  year={2013}
}

@article{cole2023quantum,
  title={On quantum optimal transport},
  author={Cole, S. and Eckstein, M. and Friedland, S. and {\.Z}yczkowski, K.},
  journal={Math. Phys. Anal. Geom.},
  volume={26},
  number={2},
  pages={14},
  year={2023}
}

@article{zhu2025unified,
  title={Unified framework for calculating convex roof resource measures},
  author={Zhu, X. and Zhang, C. and An, Z. and Zeng, B.},
  journal={npj Quantum Inf.},
  volume={11},
  number={1},
  pages={56},
  year={2025}
}

@article{beatty2025wasserstein,
  title={{W}asserstein distances on quantum structures: an overview},
  author={Beatty, E.},
  journal={arXiv preprint arXiv:2506.09794},
  year={2025}
}

@book{rockafellar1998variational,
  title={Variational analysis},
  author={Rockafellar, R. T. and Wets, R. J. B.},
  year={1998},
  publisher={Springer}
}

@book{villani2008optimal,
  title={Optimal transport: old and new},
  author={Villani, C.},
  volume={338},
  year={2008},
  publisher={Springer}
}

@article{beatty2025order,
  title={Order $p$ quantum {W}asserstein distances from couplings},
  author={Beatty, E. and Stilck Fran{\c{c}}a, D.},
  journal={Ann. Henri Poincar{\'e}},
  volume={27},
  pages={787--845},
  year={2026}
}

@article{antonini2021geometry,
  title={Geometry of {G}rassmannians and optimal transport of quantum states},
  author={Antonini, P. and Cavalletti, F.},
  journal={Ann. Sc. Norm. Super. Pisa Cl. Sci. (5)},
  volume={27},
  number={1},
  pages={201--254},
  year={2026},
  doi={10.2422/2036-2145.202302\_002}
}

@article{toth2025quantum,
  title={Quantum {W}asserstein distance and its relation to several types of fidelities},
  author={T{\'o}th, G. and Pitrik, J.},
  year={2025},
  journal={arXiv preprint arXiv:2506.14523}
}

@article{vollbrecht2001entanglement,
  title={Entanglement measures under symmetry},
  author={Vollbrecht, K. G. H. and Werner, R. F.},
  journal={Phys. Rev. A},
  volume={64},
  number={6},
  pages={062307},
  year={2001}
}

@article{choi1975completely,
  title={Completely positive linear maps on complex matrices},
  author={Choi, M.-D.},
  journal={Linear Algebra Appl.},
  volume={10},
  number={3},
  pages={285--290},
  year={1975}
}

@article{aliaga2024solution,
  title={A solution to the extreme point problem and other applications of {C}hoquet theory to {L}ipschitz-free spaces},
  author={Aliaga, R. J. and Perneck{\'a}, E. and Smith, R. J.},
  journal={arXiv preprint arXiv:2412.04312},
  year={2024}
}

@article{davidson2019noncommutative,
  title={Noncommutative {C}hoquet theory},
  author={Davidson, K. R. and Kennedy, M.},
  journal={Mem. Am. Math. Soc.},
  volume={316},
  number={1608},
  year={2025}
}

@article{brenier2020matrix,
  title={On optimal transport of matrix-valued measures},
  author={Brenier, Y. and Vorotnikov, D.},
  journal={SIAM J. Math. Anal.},
  volume={52},
  number={3},
  pages={2849--2873},
  year={2020}
}

@article{friedland2022quantum,
  title={Quantum {M}onge--{K}antorovich problem and transport distance between density matrices},
  author={Friedland, S. and Eckstein, M. and Cole, S. and {\.Z}yczkowski, K.},
  journal={Phys. Rev. Lett.},
  volume={129},
  number={11},
  pages={110402},
  year={2022}
}

@article{gozlan2017kantorovich,
  title={{K}antorovich duality for general transport costs and applications},
  author={Gozlan, N. and Roberto, C. and Samson, P.-M. and Tetali, P.},
  journal={J. Funct. Anal.},
  volume={273},
  number={11},
  pages={3327--3405},
  year={2017}
}

@article{barrett2007information,
  title={Information processing in generalized probabilistic theories},
  author={Barrett, J.},
  journal={Phys. Rev. A},
  volume={75},
  number={3},
  pages={032304},
  year={2007}
}

@article{plavala2023general,
  title={General probabilistic theories: An introduction},
  author={Pl{\'a}vala, M.},
  journal={Phys. Rep.},
  volume={1033},
  pages={1--64},
  year={2023}
}

@article{agredo2013wasserstein,
  title={A {W}asserstein-type distance to measure deviation from equilibrium of quantum {M}arkov semigroups},
  author={Agredo, J.},
  journal={Open Syst. Inf. Dyn.},
  volume={20},
  number={2},
  pages={1350009},
  year={2013}
}

@article{shirokov2006properties,
  title={Properties of probability measures on the set of quantum states and their applications},
  author={Shirokov, M. E.},
  journal={arXiv preprint arXiv:math-ph/0607019},
  year={2006}
}

@article{davidson2015choquet,
  title={The {C}hoquet boundary of an operator system},
  author={Davidson, K. R. and Kennedy, M.},
  journal={Duke Math. J.},
  volume={164},
  number={15},
  pages={2989--3004},
  year={2015}
}

@article{kennedy2022noncommutative,
  title={Noncommutative {C}hoquet simplices},
  author={Kennedy, M. and Shamovich, E.},
  journal={Math. Ann.},
  volume={382},
  pages={1591--1629},
  year={2022}
}

@article{davidson2024survey,
  title={Noncommutative {C}hoquet theory: a survey},
  author={Davidson, K. R. and Kennedy, M.},
  journal={arXiv preprint arXiv:2412.09455},
  year={2024}
}

@article{carlen2017gradient,
  title={Gradient flow and entropy inequalities for quantum {M}arkov semigroups with detailed balance},
  author={Carlen, E. A. and Maas, J.},
  journal={J. Funct. Anal.},
  volume={273},
  number={5},
  pages={1810--1869},
  year={2017}
}

@article{carlen2020noncommutative,
  title={Non-commutative calculus, optimal transport and functional inequalities in dissipative quantum systems},
  author={Carlen, E. A. and Maas, J.},
  journal={J. Stat. Phys.},
  volume={178},
  pages={319--378},
  year={2020}
}

@article{miller2026metricity,
  title={Metricity of separable quantum optimal transport with the {H}ilbert--{S}chmidt cost},
  author={Miller, T.},
  journal={arXiv preprint arXiv:2609.03769},
  year={2026}
}

@article{dusson2026wasserstein,
  title={A {W}asserstein-type metric for generic mixture models, including location-scatter and group invariant measures},
  author={Dusson, G. and Ehrlacher, V. and Nouaime, N.},
  journal={ESAIM Control Optim. Calc. Var.},
  volume={32},
  pages={19},
  year={2026},
  doi={10.1051/cocv/2026004}
}
